\documentclass{ws-ijbc}
\usepackage{ws-rotating}    % used only when sideways tables/figures are used
\usepackage{graphicx}
\usepackage{float}
\usepackage{subfigure}
\usepackage{natbib}
\usepackage{amsmath}
\usepackage{amssymb}
\usepackage{array}
\usepackage{appendix}
\usepackage{multicol}

\begin{document}

\catchline{}{}{}{}{} % Publisher's Area please ignore

\markboth{Y. et al.}{Dynamic Game of the Dual-channel Supply Chain}

\title{Dynamic Game of the Dual-channel Supply Chain under a Carbon Subsidy Policy}

\author{Yi Tian}
\address{School of Economics and Management, Tianjin University of Technology and Education\\
Tianjin, China\\
yit20160831@126.Com}

\author{Li Zhao\thanks{Email: zhao180828@126.com}} \address{Business school, Shandong Jianzhu University\\
Jinan, China\\
lzhao180828@126.com}

\author{Meihong Zhu\thanks{ Email: Mhzhu200831@126.com}} \address{Zhejiang University of Water Resources and Electric Power\\
Hang Zhou, China\\
Mhzhu200831@126.com}

\maketitle

\begin{abstract}
This study investigates the dynamic game behaviors of dual-channel supply chains involving an oligopoly manufacturer selling low-carbon products to online and offline retailers. The price game models under government subsidy are discussed under three scenarios: (1) simultaneous decision, (2) manufacturer dominates the market, and (3) retailer dominates the market. The equilibrium strategies are compared under the government subsidy policy. Using numerical simulation, complex characteristics of the dual-channel supply chain under the carbon subsidy policy are investigated. The complexity of wholesale price and sales commission of each channel is analyzed by bifurcation, largest Lyapunov exponent and basin of attraction diagrams. Furthermore, parameter adjustment and delayed feedback control methods are proven to be effective approaches to chaos control.
\end{abstract}


%\begin{multicols}{2}
\section{Introduction}
Along with the fast-evolving development of the information technique and popularization of electronic commerce, a growing number of manufacturing companies have changed their sales method from traditional retail to the deepening integration of retailing channels offline and online (\citep{tian2020coordination}). Meanwhile, the production, transportation and use of household appliances generate a large amount of carbon dioxide emission, which has attracted extensive attention from the government and society. In 2019, the retail sales of China's home appliance markets reached 8.91 trillion CNY. The retail channel of the home appliance continues to shift to the online channel, and the integration of online and offline is becoming increasingly obvious. According to the ``2019 China Home Appliance Market Report", Chinese consumers choose different channels to buy home appliances according to their preferences. Online retailing has accounted for more than 40\% of all home appliance sales, of which JD, a professional online shopping mall for low-carbon and environmental protection products in China ranks first in all channels.

In order to implement the Carbon Emission Peak of China in 2030, the government's policies have broadly supported the development green industry in recent years. China promotes the development of low-carbon energy, strengthens cooperation in the energy area, and strives to achieve this environmental sustainability commitment. The Beijing International Metropolis Clean Air Action Forum was held in 2019. At the forum, experts shared the measures, progress, medium-term and long-term plans and visions for improving air quality and reducing greenhouse gas emissions in international metropolises. They provided international experience and suggestions for Beijing to implement precise measures to effectively prevent and control air pollution and address climate change. In the latest Chinese Government Work Report, the main objectives and tasks for the ``14th Five-Year Plan" period and specific work requirements for 2021 are stated. It is proposed to accelerate the evolution of the transitional business organization and cut down on carbon emissions. Without significant technological advances and joint efforts of the whole society, it's going to be very hard to achieve carbon peak and carbon neutrality goals in China.

This paper analyzes the price competition of dual-channel supply chain in the short-term game and long-term repeated game under carbon subsidy policy. Three situations scenarios: simultaneous decision-making, manufacturer-dominated market and retailer-dominated market are discussed. And the dynamic behavior of the players is studied. The problems caused by price competition of the dual-channel supply chain under the government subsidy policy are analyzed. The numerical simulation illustrates the impact of various parameters on the price and profit of players in the supply chain.

The remainder of the research is as follows: Sec.2 provides literature related to government subsidy, price game and complexity in supply chain game. Sec.3 describes the model construction. Three dynamic game models of different power structures are studied in Sec.4. Sec.5 investigates the dynamic characteristics of the system through numerical analysis. Sec.6 analyses the system stability after chaos control. The conclusions will be given in Sec.7.

\section{Literature review}
This research draws on the dynamic price game of a dual-channel supply chain with three power structures under the carbon emission policy. Therefore, we reviewed the streams of literature on related studies in three aspects: The first stream is related to the influence of government subsidies on the supply chain game. The second stream is closely related to literature that studies the price game of the dual-channel supply chain. Finally, we study the stream which concentrates on the role of the dynamic game in a supply chain.

\subsection{Government subsidy}
The government provides low-carbon subsidies to manufacturer, retailers and consumers which can, to some extent, enhance their motivation for low-carbon production, marketing and purchase (\citep{ma2020time}).

\citep{liu2016dual} investigated the effects of prices and government subsidy on the formal channel and informal channel. It is found that government subsidy can support retailing of different channels. The product distribution of high quality products is greatly welcomed by both formal and informal channels. \citep{bian2020tax} compared net emissions, prices, profits and social welfare under two government policies (subsidy and tax) in a three-tier supply chain. \citep{li2020choice} analyzed a game theoretical approach under consumer and manufacturer subsidies. Comparing the manufacturer subsidy, the consumer subsidy brings a higher financial burden to the government and causes higher net emissions. \citep{huang2020green} discussed three scenarios in a supply chain. Under the infinite and finite loan scenarios, the optimal subsidy mode is different. \citep{yang2021technology} examined the influence of government subsidy on equilibrium solutions of a game model. \citep{wang2020duopoly} observed the influence of the marginal utility of government subsidy in the straw power supply chain. They concluded that policymakers can get a reasonable method to deal with straws.

\subsection{Price game}
Due to the rise of electronic commerce concerns, many scholars focus on studies of the diverse pricing strategies of the dual-channel supply chain(\citep{wu2021complex}). \citep{gao2020dual} presented a scenario where the government adopts an environmental subsidy to the manufacturer. The optimal prices in retail and direct channels of the supply chain were discussed based on the game theory and operations research. \citep{jabarzare2020game} studied the optimal price, quality and profit under three scenarios. \citep{ranjan2019pricing} discussed the problem of setting prices and coordination mechanism under the centralized model, decentralized model, and collaboration model of a dual-channel supply chain. \citep{matsui2020optimal} explored the optimal bargaining timing of wholesale price in a dual-channel supply chain. \citep{zhao2017pricing} considered the optimal pricing strategies of different market power structures. \citep{jamali2018game} investigated a model of pricing decision for two replacement commodities. They compared the centralized situation with the decentralized situation and concluded that the former can get a higher green degree. \citep{bayram2021order} claimed that an intuitive heuristic policy is conducive for retailers to obtain advantages in new operational challenges of the omni-channel retailing supply chain. \citep{tang2020optimising} pointed out that the retail price can be optimized by giving a new credit term-based contract. The difficulty of price games in supply chain management has intensified the conflicts and issues among channels.

\subsection{Dynamic game}
Most of the literature shows that the adjustment of dynamic price or quantity often destroys the system stability (\citep{xie2018implications}). \citep{tian2020coordination} focused on the order quantity competition in a multichannel supply chain and further discussed the complex characteristics of the system. \citep{ma2019bullwhip} investigated the influence of the bullwhip effect and complex characteristics in a dynamic price game model. \citep{mai2021dynamic} observed the impact of the system stability under decentralized decision and centralized decision in a book supply chain. \citep{wang2020duopoly} found that consumers' green preference directly affects social stability. At present, most manufacturers sell their products through dual-channel. The complexity of dual-channel supply chain has attracted widespread academic attention. It was demonstrated that the prices, quantities, and profits were significantly affected by the stable performance of the system. Therefore, a large number of scholars gravitate to using nonlinear dynamical methods in the analysis of practical problems (\citep{ma2022evolution}; \citep{ma2022optimal}; \citep{chen2020q}; \citep{wu2018equilibrium}; \citep{chen2016dynamics}).

Although there exist some studies about the dual-channel supply chain in the literature, few of them focus on the analysis about the impact of the carbon subsidy policy and different power structures. Furthermore,  few academic research on the influence of pricing system stability. If the lack of system stability, the players cannot obtain stable and long-run profits. From an application perspective, the system stability helps the players obtain long-run profits. The relevant literature is summarized in Table 1.
%±í¸ñÒ»
\begin{table}[ht]
\centering
\tbl{Control effect comparison.}
{\begin{tabular}{ c c c c c c c }
\hline
 Study article & Dual-channel & Pricing decision & Nash &  Stackelberg & Government
Regulation & Dynamic game\\
\hline
\citep{mai2021dynamic} & $\times$& $\surd$ &$\times$&$\surd$&$\times$ &$ \surd$ \\
\citep{bian2020tax}  & $ \times $ & $\surd$ & $\times$  &$\surd$& $\surd$  & $\times$ \\
\citep{li2020choice} &  $\times$  & $\surd$ & $\times$  &$\surd$& $\surd$  & $\times$\\
\citep{huang2020green} & $\times$ & $\surd$ & $\times$  &$\surd$& $\surd$  & $\times$\\
\citep{gao2020dual} & $\surd$  &$ \surd $& $\times $ &$\surd$& $\surd$  & $\times$\\
\citep{jabarzare2020game}  &$ \surd$  & $\surd$ & $\times$  &$\surd$& $\times$  & $\times$\\
\citep{ranjan2019pricing} &$ \surd$  & $\surd$ & $\times$  &$\surd$& $\times$  & $\times$\\
\citep{matsui2020optimal}  &$ \surd$  & $\surd$ & $\surd$  &$\surd$& $\times$  & $\times$\\
\citep{wang2020duopoly} &$\times$   & $\surd$ & $\times$   &$\surd$& $\times$  & $\times$\\
\citep{jamali2018game} &$\surd$  & $\surd$ & $\surd$   &$\surd$& $\times$  & $\times$\\
\citep{zhao2017pricing} &$\surd$  & $\surd$ & $\surd$   &$\surd$& $\times$  & $\times$\\
This paper &$\surd$   & $\surd$ & $\surd$   &$\surd$& $\surd$  & $\surd$\\
\hline
\end{tabular}}
\end{table}

Therefore, we summarize our main contributions as follows: Firstly, a price game model of dual-channel supply chain model incorporating the carbon subsidy policy is proposed. Secondly, the complexity of the dual-channel supply chain under three dynamic game structures is investigated. Finally, it provides valuable suggestions for the managers of enterprises in the dynamically adjusting strategies of wholesale price and sales commission under the three power structures. Table 2 shows the parameters in the paper.
%±í¸ñ¶þ
\begin{table}[h]
\tbl{Parameters.}
{\begin{tabular}{p{2.0in}p{2.9in}} \hline
Parameter\&Superscript &Definition \\ \hline
\qquad\qquad  { $ q_1 $ } & Offline retailer's sales quantity \\
\qquad\qquad  { $ q_2 $ } & Online retailer's sales quantity \\
\qquad\qquad  { $ p_1 $ } & Offline retailer's sales price \\
\qquad\qquad  { $ p_2 $ } & Online retailer's sales price \\
\qquad\qquad  { $ c $ } & Production costs \\
\qquad\qquad  { $ \mu $ } & Two retailers' substitution coefficient \\
\qquad\qquad  { $ \eta $ } & Consumers' channel preference \\
\qquad\qquad  { $ \lambda $ } & Consumers' low-carbon preference \\
\qquad\qquad  { $ g $ } & Manufacturer's low-carbon innovativeness \\
\qquad\qquad  { $ g_0 $ } & Threshold for government subsidy for low-carbon products \\
\qquad\qquad  { $ l $ } & Level of government low-carbon subsidy to manufacturer \\
\qquad\qquad  { $ \theta $ } & The proportion of offline retailer undertaking low-carbon innovation \\
\qquad\qquad  { $ g_1 $ } & Adjustment speed of offline retailer's sales commissions \\
\qquad\qquad  { $ g_2 $ } & Adjustment speed of online retailer's sales commissions \\
\qquad\qquad  { $ g_3 $ } & Adjustment speed of manufacturer's wholesale prices with regard to offline channel \\
\qquad\qquad  { $ g_4 $ } & Adjustment speed of manufacturer's wholesale prices with regard to online channel\\
\qquad\qquad  { $ \pi_m $ } & Manufacturer's profits \\
\qquad\qquad  { $ \pi_{r1} $ } & Offline retailer's profits \\
\qquad\qquad  { $ \pi_{r2} $ } & Online retailer's profits \\
\qquad\qquad  { $ w_1 $ } & Wholesale prices for offline retailer \\
\qquad\qquad  { $ w_2 $ } & Wholesale prices for online retailer \\
\qquad\qquad  { $ k_1 $ } & Offline retailer's sales commissions \\
\qquad\qquad  { $ k_2 $ } & Online retailer's sales commissions \\
\qquad\qquad  { $ NG $ } & Nash game \\
\qquad\qquad  { $ MS $ } & Manufacturer-dominated Stackelberg game \\
\qquad\qquad  { $ RS $ } & Retailer-dominated Stackelberg game \\\hline

\end{tabular}}
\end{table}

\section{Model Construction}\label{sec:basicmodel}

The dual-channel supply chain system consists of a manufacturer, an offline retailer (Retailer 1), and an online retailer (Retailer 2). The manufacturer provides green products to retailer 1 and retailer 2 at the same unit cost ($c$) and wholesale price($w_i$). Subsequently, the two retailers sell their merchandise through online channel and traditional channel respectively. Retailers determine their unit sales commission $k_i(i=1,2)$. The retail price is $p_i=w_i+k_i$. In addition, the competitive pressure faced by online and offline retailers is increasing. And they intensely compete in the market shared by both.

The replaceability of two retailers in the market is expressed by a parameter $\mu$. $\mu=0$ represents the monopoly of one retailer in the market, and $\mu=1$ represents the existence of two completely replaceable retailers. The demand function can be shown as follows:

\begin{equation}
q_1=1-p_1+\mu{p_2}, \\
q_2=\eta-p_2+\mu{p_1}\\
\end{equation}

Next, we will study the price game under government low-carbon subsidies. Due to consumers prefer green environmental protection commodities, the order quantity of low-carbon products will accordingly increase. The demand function can be obtained as follows:

\begin{equation}
q_1=1-p_1+\mu{p_2}+\lambda g, \\
q_2=\eta-p_2+\mu{p_1}+\lambda g\\
\end{equation}

Where, $\lambda$ is the consumers' low-carbon preference, and $g$ is the manufacturer's low-carbon innovativeness.

The wholesale prices are $w_1$ and $w_2$, while the sales commissions are $k_1$ and $k_2$. The manufacturer can obtain the government subsidy when they engage in low-carbon innovation activities. The threshold of government subsidy for low-carbon products is $g_0$. The government subsidy can be obtained only when the low-carbon innovativeness of products is greater than $g_0$. Therefore, the amount of government subsidy is a quadratic function $\ell(g-g_0)^2$, where $\ell$ is the level of government low-carbon subsidy to the manufacturer. The marginal utility of the manufacturer's innovative low-carbon investment gradually decreases, and the amount of the manufacturer's investment in low-carbon technology is $\frac{1}{2}\lambda g^2$. Offline retailer strives to bear the corporation social responsibility and takes the initiative to share manufacturer's low-carbon innovation costs. The proportion of offline retailer undertaking low-carbon innovation is $\theta$.

The operational profits of the manufacturer and the two retailers are shown in Eq.(3), Eq.(4) and Eq.(5).
\begin{equation}
 \begin{array}{l}
\pi_{r 1}=k_{1}\left(1-\left(k_{1}+w_{1}\right)+\mu\left(k_{2}+w_{2}\right)+\lambda g\right)-\frac{1}{2} \theta \lambda g^{2} \\
\end{array}
\end{equation}
\begin{equation}
\pi_{r 2}=k_{2}\left(\eta-\left(k_{2}+w_{2}\right)+\mu\left(k_{1}+w_{1}\right)+\lambda g\right) \\
\end{equation}
\begin{equation}
\pi_{m}=w_{1}\left(1-\left(k_{1}+w_{1}\right)+\mu\left(k_{2}+w_{2}\right)+\lambda g\right)+w_{2}\left(\eta-\left(k_{2}+w_{2}\right)+\mu\left(k_{1}+w_{1}\right)+\lambda g\right)-\frac{1}{2}(1-\theta) \lambda g^{2}+\ell\left(g-g_{0}\right)^{2}
\end{equation}

Considering the different power structures of manufacturer and retailers in the decision-making process, three cases are discussed below: the Nash game in which manufacturer and retailers have equal power structures; the manufacturer-dominated Stackelberg game which means the manufacturer is dominant and two retailers are subordinate; the retailer-dominated Stackelberg game in which two retailers are dominant in the market and manufacturer is a follower.

\subsection{Nash game model}

Firstly, we will consider the situation that manufacturer and retailers have equal power structures. Players can make decisions simultaneously. In this paper, the Nash game model is referred to as model NG and denoted by the superscript N. When manufacturer and online and offline retailers are equally matched, they simultaneously decide wholesale prices $w_1$ and $w_2$ and sales commissions $k_1$ and $k_2$ in order to maximize profits.

The first-order partial derivatives of equations (3), (4) and (5) can be calculated:

\begin{equation}
\left\{\begin{array}{l}
\frac{\partial \pi_{r 1}}{\partial k_{1}}=1-2 k_{1}-w_{1}+\lambda g+\left(k_{2}+w_{2}\right) \mu \\
\frac{\partial \pi_{r 2}}{\partial k_{2}}=-2 k_{2}-w_{2}+\lambda g+\eta+\left(k_{1}+w_{1}\right) \mu \\
\frac{\partial \pi_{m}}{\partial w_{1}}=1-k_{1}-2 w_{1}+\lambda g+w_{2} \mu+\left(k_{2}+w_{2}\right) \mu \\
\frac{\partial \pi_{m}}{\partial w_{2}}=-k_{2}-2 w_{2}+\lambda g+\eta+w_{1} \mu+\left(k_{1}+w_{1}\right) \mu
\end{array}\right.
\end{equation}

By making the first-order partial derivatives of manufacturer and retailers zero simultaneous, the equilibrium solution of the decision is obtained as follows:

\begin{equation}
\left\{\begin{array}{l}
k_{1}^{N^{*}}=-\frac{3+\eta \mu+g \lambda(3+\mu)}{-9+\mu^{2}} \\
k_{2}^{N^{*}}=-\frac{3 \eta+\mu+g \lambda(3+\mu)}{-9+\mu^{2}} \\
w_{1}^{N^{*}}=\frac{3+4 \eta \mu+\mu^{2}+g \lambda\left(3+4 \mu+\mu^{2}\right)}{9-10 \mu^{2}+\mu^{4}} \\
w_{2}^{N^{*}}=\frac{4 \mu+\eta\left(3+\mu^{2}\right)+g \lambda\left(3+4 \mu+\mu^{2}\right)}{9-10 \mu^{2}+\mu^{4}}
\end{array}\right.
\end{equation}

\subsection{Manufacturer-dominated game model}

Then we consider the situation in which the manufacturer is the leading position in the market. The manufacturer is the first mover who decides the wholesale price by maximizing the profit. And the two retailers are in a subordinate position. The retailers determine the sales commissions by maximizing their profit according to the leader's decision. This model is referred to as model MS, which is labelled by superscript M.
According to the reverse order solution method, the first-order partial derivative of the two retailers' profits is zero to obtain simultaneously:

\begin{equation}
\left\{\begin{array}{l}
k_{1}^{M^{*}}=-\frac{2+w_{2} \mu+\eta \mu+g \lambda(2+\mu)+w_{1}\left(-2+\mu^{2}\right)}{-4+\mu^{2}} \\
k_{2}^{M^{*}}=-\frac{2 \eta+2 g \lambda+\mu+w_{1} \mu+g \lambda \mu+w_{2}\left(-2+\mu^{2}\right)}{-4+\mu^{2}}
\end{array}\right.
\end{equation}

We bring Eq. (8) into Eq. (5) and then obtain the first partial derivatives of the manufacturer's profit.

\begin{equation}
\left\{\begin{array}{l}
\frac{\partial \pi_{m}}{\partial w_{1}}=-\frac{2+2 w_{2} \mu+\eta \mu+g \lambda(2+\mu)+2 w_{1}\left(-2+\mu^{2}\right)}{-4+\mu^{2}} \\
\frac{\partial \pi_{m}}{\partial w_{2}}=-k_{2}-2 w_{2}+\lambda g+\eta+w_{1} \mu+\left(k_{1}+w_{1}\right) \mu
\end{array}\right.
\end{equation}

$w_1^{M*}$and $w_2^{M*}$ can be solved by making Eq.(9) to be zero. Combining Eq.(8), the expression of the equilibrium solution is as follows:

\begin{equation}
\left\{\begin{array}{l}
k_{1}^{M^{*}}=\frac{2+\eta \mu+g \lambda(2+\mu)}{8-2 \mu^{2}} \\
k_{2}^{M^{*}}=\frac{2 \eta+\mu+g \lambda(2+\mu)}{8-2 \mu^{2}} \\
w_{1}^{M^{*}}=\frac{1+\eta \mu+g \lambda(1+\mu)}{2-2 \mu^{2}} \\
w_{2}^{M^{*}}=\frac{\eta+\mu+g \lambda(1+\mu)}{2-2 \mu^{2}}
\end{array}\right.
\end{equation}

\subsection{Retailer-dominated game model}
We further consider the situation in which the two retailers have an advantageous position in price competition and use the dominant position of the supply chain to decide their respective sales commissions simultaneously. The manufacturer decides wholesale price according to the decision of the retailer. The retailer-dominated Stackelberg game model is abbreviated as model RS, denoted by superscript R.
The first partial derivative of the manufacturer's profit is zero, and then the following equations are obtained:

\begin{equation}
\left\{\begin{array}{l}
w_{1}^{R^{*}}=\frac{1+\eta \mu+g \lambda(1+\mu)+k_{1}\left(-1+\mu^{2}\right)}{2-2 \mu^{2}} \\
w_{2}^{R^{*}}=-\frac{\eta+g \lambda+\mu+g \lambda \mu+k_{2}\left(-1+\mu^{2}\right)}{2\left(-1+\mu^{2}\right)}
\end{array}\right.
\end{equation}

Incorporating equation (11) into equations (3) and (4) to find the first derivative of the two retailers' profits:

\begin{equation}
\left\{\begin{array}{l}
\frac{\partial \pi_{r 1}}{\partial k_{1}}=\frac{1}{2}\left(1-2 k_{1}+g \lambda+k_{2} \mu\right) \\
\frac{\partial \pi_{r 2}}{\partial k_{2}}=\frac{1}{2}\left(-2 k_{2}+\eta+g \lambda+k_{1} \mu\right)
\end{array}\right.
\end{equation}

$k_1^{R*}$ and $k_2^{R*}$ can be solved by making the first order derivative of Eq.(12) to be zero. Combining Eq.(11), $w_1^{R*}$and $w_2^{R*}$ can be obtained. The expressions for $k_1^{R*}$, $k_2^{R*}$, $w_1^{R*}$ and $w_2^{R*}$ are as follows:

\begin{equation}
\left\{\begin{array}{l}
k_{1}^{R^{*}}=-\frac{2+\eta \mu+g \lambda(2+\mu)}{-4+\mu^{2}} \\
k_{2}^{R^{*}}=-\frac{2 \eta+\mu+g \lambda(2+\mu)}{-4+\mu^{2}} \\
w_{1}^{R^{*}}=\frac{2+3 \eta \mu+\mu^{2}+g \lambda\left(2+3 \mu+\mu^{2}\right)}{2\left(4-5 \mu^{2}+\mu^{4}\right)} \\
w_{2}^{R^{*}}=\frac{3 \mu+\eta\left(2+\mu^{2}\right)+g \lambda\left(2+3 \mu+\mu^{2}\right)}{2\left(4-5 \mu^{2}+\mu^{4}\right)}
\end{array}\right.
\end{equation}

\subsection{Comparative analysis of single stage game}
Next, the single-stage game of the three models is analyzed, and the influence of enterprise equilibrium strategy on enterprise decision-making is compared.

%±í¸ñÈý
\begin{table}[ht]
\centering
\tbl{Variation of equilibrium strategy with parameters.}
{\begin{tabular}{ c c c c c c c }
\hline
 Model & Equilibrium solution & g & $\lambda$ & $\eta$  \\
\hline
Model NG & $ k_{1}^{N^{*}}$ & $ + $& $ + $ &$ + $ \\
         & $ k_{1}^{N^{*}}$ & $ + $& $ + $ &$ + $ \\
         & $ w_{1}^{N^{*}}$ & $ + $& $ + $ &$ + $ \\
         & $ w_{1}^{N^{*}}$ & $ + $& $ + $ &$ + $ \\
\hline
Model MS  & $ k_{1}^{M^{*}}$ & $ + $& $ + $ &$ + $ \\
          & $ k_{1}^{M^{*}}$ & $ + $& $ + $ &$ + $ \\
          & $ w_{1}^{M^{*}}$ & $ + $& $ + $ &$ + $ \\
          & $ w_{1}^{M^{*}}$ & $ + $& $ + $ &$ + $ \\
\hline
Model RS  & $ k_{1}^{R^{*}}$ & $ + $& $ + $ &$ + $ \\
          & $ k_{1}^{R^{*}}$ & $ + $& $ + $ &$ + $ \\
          & $ w_{1}^{R^{*}}$ & $ + $& $ + $ &$ + $ \\
          & $ w_{1}^{R^{*}}$ & $ + $& $ + $ &$ + $ \\
\hline
\end{tabular}}
\end{table}

Table 3 indicates the changing situation of equilibrium strategy values with the different parameters under the three models. We obtain the equilibrium solution and analyze the variation of the equilibrium solution with the parameters $g$, $\lambda$ and $\eta$. It can be seen from the table 3 that no matter which model is used, the equilibrium solutions of wholesale price and sales commission are positively related to the degree of green products $g$, consumers' low-carbon preference $\lambda$, and channel preference $\eta$. It can be seen that if low-carbon innovativeness of products is relatively high and consumers pay more attention to and prefer more green low-carbon products, manufacturer and retailer can moderately increase product prices. Consumers' preference for the online channel can also make the supply chain players moderately increase product prices.
%±í¸ñËÄ
\begin{table}[ht]
\centering
\tbl{Equilibrium strategy for three models.}
{\begin{tabular}{ c c c c c c c c c c }
\hline
Model & $ k_{1}^{*} $ & $ k_{2}^{*}$ & $ w_{1}^{*}$ & $ w_{2}^{*}$ & $ \pi_{r1}^{*}$ & $ \pi_{r2}^{*}$ & $ \pi_{m}^{*}$ & $ \pi_{SC}^{*}$ \\
\hline
Model NG &  $ 0.426 $  &  $ 0.454 $  &  $ 0.87 $  &  $ 0.89 $  &   $0.179 $  & $ 0.206 $  &  $ 0.774  $ &  $ 1.159  $ \\
Model MS &  $ 0.357  $ &  $ 0.377 $  &  $ 1.083 $ & $  1.117  $ &  $ 0.125  $ &  $ 0.142  $ &  $ 0.806 $  & $ 1.073 $  \\
Model RS &  $ 0.713  $ &  $ 0.753 $  &  $ 0.727 $ &  $  0.74  $ &  $ 0.252  $ &  $ 0.284  $ &  $ 0.537  $ &  $ 1.073 $  \\
\hline
\end{tabular}}
\end{table}

Based on the proposition, we can set the parameters as $\mu=0.5$, $\eta=1.1$, $\lambda=0.1$, $g=0.5$, $g_0=0.1$ and $l=0.1$. Table 4 shows the comparison of enterprise equilibrium strategies under Model NG, Model MS and Model RS. When the government subsidizes low-carbon innovative enterprise and offline retailer share the cost of low-carbon innovation, the situation of the price system and profit of the supply chain system can be seen in Table 4. Whoever dominates the market will hold the initiative of market pricing, and the pricing of this party will continue to increase appropriately. The profit of this party will also increase more appropriately than in the other two cases. Regarding the total profit, it is the highest when the players make decisions simultaneously than in the other two situations.

\section{Research on dynamic game model}
In this section, based on government subsidy, we will study the complexity of the dynamic game process for this model with different power structures. The real market is increasingly complex and changeable. Players in the supply chain cannot completely obtain market information. Therefore, market decision-makers generally adopt bounded rational expectations, and the decision-making $t$ period to $t+1$ period is shown in equation (14):

\begin{equation}
\left\{\begin{array}{l}
k_{1, t+1}=k_{1, t}+g_{1} k_{1, t} \frac{\partial E\left(\pi_{r 1,t}\right)}{\partial k_{1, t}} \\
k_{2, t+1}=k_{2, t}+g_{2} k_{2, t} \frac{\partial E\left(\pi_{r 2,t}\right)}{\partial k_{2, t}} \\
w_{1, t+1}=w_{1, t}+g_{3} w_{1, t} \frac{\partial E\left(\pi_{m, t}\right)}{\partial w_{1, t}} \\
w_{2, t+1}=w_{2, t}+g_{4} w_{2, t} \frac{\partial E\left(\pi_{m, t}\right)}{\partial w_{2, t}}
\end{array}\right.
\end{equation}

The price adjustment speed includes two aspects: the adjustment speeds of wholesale price and sales commission. In equation (14), $g_{1}$ and $g_{2}$ are the adjustment speeds of the two retailers' sales commissions, respectively.  $g_{3}$ and $g_{4}$ are the adjustment speeds of the manufacturer's wholesale prices. $g_{i}(i=1,2,3 $ and $ 4)$ is larger than zero. Both manufacturer and retailers make adjustments of bounded rational price based on the current marginal profit to decide whether to increase or cut the price in the next period. The dynamic adjustment process of the three models is discussed below.
\subsection{Nash dynamic game model}
The repeated game model for supply chain members to make simultaneous decisions based on bounded rational expectations is as follows:

\begin{equation}
\left\{\begin{array}{l}
k_{1, t+1}=k_{1, t}+g_{1} k_{1, t}\left(1-2 k_{1, t}-w_{1, t}+g \lambda+\left(k_{2, t}+w_{2, t}\right) \mu\right) \\
k_{2, t+1}=k_{2, t}+g_{2} k_{2, t}\left(-2 k_{2, t}-w_{2, t}+\eta+g \lambda+\left(k_{1, t}+w_{1, t}\right) \mu\right) \\
w_{1, t+1}=w_{1, t}+g_{3} w_{1, t}\left(1-k_{1, t}-2 w_{1, t}+g \lambda+k_{2, t} \mu+2 w_{2, t} \mu\right) \\
w_{2, t+1}=w_{2, t}+g_{4} w_{2, t}\left(-k_{2, t}-2 w_{2, t}+\eta+g \lambda+k_{1, t} \mu+2 w_{1, t} \mu\right)
\end{array}\right.
\end{equation}

Eq.(15) is set to be $k_{i, t+1}=k_{i, t}$,$\quad w_{i, t+1}=w _{i, t}$. We can calculate 16 equilibrium points. $E^{N^{*}}\left(k_{1}^{N^{*}} k_{2}^{N^{*}} w_{1}^{N^{*}} w_{2}^{N^{*}}\right)$ is the Nash equilibrium point. The Jacobian matrix is given as follow.

\begin{equation}
J^{N}=\left(\begin{array}{cccc}
J_{11} & g_{1} k_{1} \mu & -g_{1} k_{1} \mu & g_{1} k_{1} \mu \\
g_{2} k_{2} \mu & J_{22} & g_{2} k_{2} \mu & -g_{2} k_{2} \\
-g_{3} w_{1} & g_{3} w_{1} \mu & J_{33} & 2 g_{3} w_{1} \mu \\
g_{4} w_{2} \mu & -g_{4} w_{2} & 2 g_{4} w_{2} \mu & J_{44}
\end{array}\right)
\end{equation}

where,

\begin{equation}
\begin{array}{l}
J_{11}^{N}=1-2 g_{1} k_{1}+g_{1}\left(1-2 k_{1}-w_{1}+g \lambda+\left(k_{2}+w_{2}\right) \mu\right) \\
J_{22}^{N}=1-2 g_{2} k_{2}+g_{2}\left(-2 k_{2}-w_{2}+g \lambda+\eta+\left(k_{1}+w_{1}\right) \mu\right) \\
J_{33}^{N}=1-2 g_{3} w_{1}+g_{3}\left(1-k_{1}-2 w_{1}+g \lambda+k_{2} \mu+2 w_{2} \mu\right) \\
J_{44}^{N}=1-2 g_{4} w_{2}+g_{4}\left(-k_{2}-2 w_{2}+\eta+g \lambda+k_{1} \mu+2 w_{1} \mu\right)
\end{array}
\end{equation}

\subsection{Manufacturer-dominated Stackelberg dynamic game model}

In the manufacturer-dominated Stackelberg game, the two retailers get the marginal profit with respect to $k_1$ and $k_2$ at $t$ period; In the $t+1$ period, retailers' reaction functions $k_1^*$ and $k_2^*$ can be substituted into equation (5) to obtain the marginal profit equation. The manufacturer makes the bounded rationality expectation decision based on equation (9) to obtain the following equation.

\begin{equation}
\left\{\begin{array}{l}
w_{1, t+1}=w_{1, t}+g_{3} w_{1, t}\left(-\frac{2+2 w_{2, t} \mu+\eta \mu+g \lambda(2+\mu)+2 w_{1, t}\left(-2+\mu^{2}\right)}{-4+\mu^{2}}\right) \\
w_{2, t+1}=w_{2, t}+g_{4} w_{2, t}\left(-k_{2, t}-2 w_{2, t}+g \lambda+\eta+w_{1, t} \mu+\left(k_{1, t}+w_{1, t}\right) \mu\right)
\end{array}\right.
\end{equation}

After knowing $w_{1,t+1}$ and $w_{2,t+1}$, the two retailers make decisions $k_{1,t+1}$ and $k_{2,t+1}$. By combining with (18), manufacturer-dominated Stackelberg dynamic game system can be obtained:

\begin{equation}
\left\{\begin{array}{l}
w_{1, i+1}=w_{1, t}+g_{3} w_{1, t}\left(-\frac{2+2 w_{2, t} \mu+\eta \mu+g \lambda(2+\mu)+2 w_{1, t}\left(-2+\mu^{2}\right)}{-4+\mu^{2}}\right) \\
w_{2, t+1}=w_{2, t}+g_{4} w_{2, t}\left(-k_{2, t}-2 w_{2, t}+g \lambda+\eta+w_{1, t} \mu+\left(k_{1, t}+w_{1, t}\right) \mu\right) \\
k_{1, t+1}=k_{1, t}+g_{1} k_{1, t}\left(1-2 k_{1, t}-w_{1, t+1}+g \lambda+\left(k_{2, t}+w_{2, t+1}\right) \mu\right) \\
k_{2, t+1}=k_{2, t}+g_{2} k_{2, t}\left(-2 k_{2, t}-w_{2, t+1}+\eta+g \lambda+\left(k_{1, t}+w_{1, t+1}\right) \mu\right)
\end{array}\right.
\end{equation}

When $k_{i,t+1}=k_{i,t}$,$w_{i,t+1}=w_{i,t}$, the Nash equilibrium point    $E^{M^{*}}\left(k_{1}^{M^{*}} k_{2}^{M^{*}} w_{1}^{M^{*}} w_{2}^{M^{*}}\right)$can be obtained.

The Jacobi matrix can be given as:

\begin{equation}
J^{M}=\left(\begin{array}{lllc}
J_{11} & g_{1} k_{1} \mu & J_{13} & J_{14} \\
g_{2} k_{2} \mu & J_{22} & J_{23} & J_{24} \\
0 & 0 & J_{33} & -\frac{2 g_{3} w_{1} \mu}{-4+\mu^{2}} \\
0 & 0 & -\frac{2 g_{4} w_{2} \mu}{-4+\mu^{2}} & J_{44}
\end{array}\right)
\end{equation}

where,

\begin{equation}
\begin{array}{l}
J_{11}^{M}=1+g_{1}-4 g_{1} k_{1}-g_{1} w_{1}+g g_{1} \lambda+\frac{g_{1} g_{3} w_{1}\left(2+2 w_{2} \mu+\eta \mu+g \lambda(2+\mu)+2 w_{1}\left(-2+\mu^{2}\right)\right)}{-4+\mu^{2}}\\
+g_{1} \mu\left(k_{2}+w_{2}-\frac{g_{4} w_{2}\left(2 \eta+2 g \lambda+\mu+2 w_{1} \mu+g \lambda \mu+2 w_{2}\left(-2+\mu^{2}\right)\right)}{-4+\mu^{2}}\right)\\
J_{13}^{M}=\frac{g_{1} k_{1}\left(4-\left(1+2 g_{4} w_{2}\right) \mu^{2}+g_{3}\left(2+2 w_{2} \mu+\eta \mu+g \lambda(2+\mu)+4 w_{1}\left(-2+\mu^{2}\right)\right)\right)}{-4+\mu^{2}}\\
J_{14}^{M}=-\frac{g_{1} k_{1} \mu\left(4-2 g_{3} w_{1}-\mu^{2}+g_{4}\left(2 \eta+2 g \lambda+\mu+2 w_{1} \mu+g \lambda \mu+4 w_{2}\left(-2+\mu^{2}\right)\right)\right)}{-4+\mu^{2}}\\
J_{22}^{M}=1-4 g_{2} k_{2}-g_{2} w_{2}+g_{2} \eta+g g_{2} \lambda+\frac{g_{2} g_{4} w_{2}\left(2 \eta+2 g \lambda+\mu+2 w_{1} \mu+g \lambda \mu+2 w_{2}\left(-2+\mu^{2}\right)\right)}{-4+\mu^{2}}\\
+g_{2} \mu\left(k_{1}+w_{1}-\frac{g_{3} w_{1}\left(2+2 w_{2} \mu+\eta \mu+g \lambda(2+\mu)+2 w_{1}\left(-2+\mu^{2}\right)\right)}{-4+\mu^{2}}\right)\\
J_{23}^{M}=-\frac{g_{2} k_{2} \mu\left(4-2 g_{4} w_{2}-\mu^{2}+g_{3}\left(2+2 w_{2} \mu+\eta \mu+g \lambda(2+\mu)+4 w_{1}\left(-2+\mu^{2}\right)\right)\right)}{-4+\mu^{2}}\\
J_{24}^{M}=\frac{g_{2} k_{2}\left(4-\left(1+2 g_{3} w_{1}\right) \mu^{2}+g_{4}\left(2 \eta+2 g \lambda+\mu+2 w_{1} \mu+g \lambda \mu+4 w_{2}\left(-2+\mu^{2}\right)\right)\right)}{-4+\mu^{2}}\\
J_{33}^{M}=\frac{-4+\mu^{2}-g_{3}\left(2+2 w_{2} \mu+\eta \mu+g \lambda(2+\mu)+4 w_{1}\left(-2+\mu^{2}\right)\right)}{-4+\mu^{2}}\\
J_{44}^{M}=\frac{-4+\mu^{2}-g_{4}\left(2 \eta+2 g \lambda+\mu+2 w_{1} \mu+g \lambda \mu+4 w_{2}\left(-2+\mu^{2}\right)\right)}{-4+\mu^{2}}
\end{array}
\end{equation}

\subsection{Retailer-dominated Stackelberg dynamic game model}

For the retailer-dominated Stackelberg dynamic game, in the $t$ period, the manufacturer obtains $w_1$ and $w_2$ and substitutes $w_1^*$ and $w_2^*$ into the profit formulas (3) and (4) to obtain marginal profit (12). The retailers make a bounded rational expected decision based on formula (12) to obtain formula (22).

\begin{equation}
\left\{\begin{array}{l}
k_{1, t+1}=k_{1, t}+g_{1} k_{1, t}\left(\frac{1}{2}\left(1-2 k_{1, t}+g \lambda+k_{2, t} \mu\right)\right) \\
k_{2, t+1}=k_{2, t}+g_{2} k_{2, t}\left(\frac{1}{2}\left(-2 k_{2, t}+\eta+g \lambda+k_{1, t} \mu\right)\right)
\end{array}\right.
\end{equation}

The manufacturer makes decisions according to bounded rational expectations based on knowing the decision of retailer 1 and retailer 2. Combining with (22), the retailer-dominated Stackelberg dynamic game system can be obtained as:

\begin{equation}
\left\{\begin{array}{l}
k_{1, t+1}=k_{1, t}+g_{1} k_{1, t}\left(\frac{1}{2}\left(1-2 k_{1, t}+g \lambda+k_{2, t} \mu\right)\right) \\
k_{2, t+1}=k_{2, t}+g_{2} k_{2, t}\left(\frac{1}{2}\left(-2 k_{2, t}+\eta+g \lambda+k_{1, t} \mu\right)\right) \\
w_{1, t+1}=w_{1, t}+g_{3} w_{1, t}\left(1-k_{1, t+1}-2 w_{1, t}+g \lambda+k_{2, t+1} \mu+2 w_{2, t} \mu\right) \\
w_{2, t+1}=w_{2, t}+g_{4} w_{2, t}\left(-k_{2, t+1}-2 w_{2, t}+\eta+g \lambda+k_{1, t+1} \mu+2 w_{1, t} \mu\right)
\end{array}\right.
\end{equation}

The Nash equilibrium point $E^{R^{*}}\left(k_{1}^{R^{*}} k_{2}^{R^{*}} w_{1}^{R^{*}} w_{2}^{R^{*}}\right)$ can be obtained as $k_{i,t+1}=k_{i,t}$ and $w_{i,t+1}=w_{i,t}$.
The Jacobian matrix of  Eq. (26) can be written as:

\begin{equation}
J^{R}=\left(\begin{array}{cccc}
J_{11} & \frac{g_{1} k_{1} \mu}{2} & 0 & 0 \\
g_{2} k_{2} \mu & J_{22} & 0 & 0 \\
2 & J_{32} & J_{33} & 2 g 3 w 1 \mu \\
J_{31} & J_{42} & 2 g_{4} w_{2} \mu & J_{44}
\end{array}\right)
\end{equation}

where,

\begin{equation}
\begin{array}{l}
J_{11}^{R}=1+\frac{1}{2} g_{1}\left(1-4 k_{1}+g \lambda+k_{2} \mu\right) \\
J_{22}^{R}=1+\frac{1}{2} g_{2}\left(-4 k_{2}+\eta+g \lambda+k_{1} \mu\right) \\
J_{31}^{R}=\frac{1}{2} g_{3} w_{1}\left(-2+g_{2} k_{2} \mu^{2}+g_{1}\left(-1+4 k_{1}-g \lambda-k_{2} \mu\right)\right) \\
J_{32}^{R}=\frac{1}{2} g_{3} w_{1} \mu\left(2-g_{1} k_{1}+g_{2}\left(-4 k_{2}+\eta+g \lambda+k_{1} \mu\right)\right) \\
J_{33}^{R}=1+g_{3}-g_{3} k_{1}-4 g_{3} w_{1}+g g_{3} \lambda+2 g_{3} w_{2} \mu-\frac{1}{2} g_{1} g_{3} k_{1}\left(1-2 k_{1}+g \lambda+k_{2} \mu\right) \\
+g_{3} \mu\left(k_{2}+\frac{1}{2} g_{2} k_{2}\left(-2 k_{2}+\eta+g \lambda+k_{1} \mu\right)\right) \\
J_{41}^{R}=\frac{1}{2} g_{4} w_{2} \mu\left(2-g_{2} k_{2}+g_{1}\left(1-4 k_{1}+g \lambda+k_{2} \mu\right)\right) \\
J_{42}^{R}=\frac{1}{2} g_{4} w_{2}\left(-2+g_{1} k_{1} \mu^{2}+g_{2}\left(4 k_{2}-\eta-g \lambda-k_{1} \mu\right)\right) \\
J_{44}^{R}=1-2 g_{4} w_{2}+g_{4}\left(-k_{2}-2 w_{2}+\eta+g \lambda+2 w_{1} \mu-\frac{1}{2} g_{2} k_{2}\left(-2 k_{2}+\eta+g \lambda+k_{1} \mu\right)\right) \\
+\mu\left(k_{1}+\frac{1}{2} g_{1} k_{1}\left(1-2 k_{1}+g \lambda+k_{2} \mu\right)\right)
\end{array}
\end{equation}

The characteristic polynomial is:
\begin{equation}
F(A)=A^{4}+\zeta_{3} A^{3}+\zeta_{2} A^{2}+\zeta_{1} A+\zeta_{0}
\end{equation}

To keep the stable region of the dynamic system, the following Eq.(27) should be met:

\begin{equation}
\left\{\begin{array}{l}
1+\zeta_{3}+\zeta_{2}+\zeta_{1}+\zeta_{0}>0 \\
1-\zeta_{3}+\zeta_{2}-\zeta_{1}+\zeta_{0}<0 \\
\left|\zeta_{0}\right|<1 \\
\left|\kappa_{0}\right|<\left|\kappa_{3}\right| \\
\left|\chi_{0}\right|<\left|\chi_{2}\right|
\end{array}\right.
\end{equation}

The detailed expressions of $\kappa_0$,$\kappa_3$,$\chi_0$ and $\chi_2$ are as follows:
\begin{equation}
\begin{array}{l}
\kappa_{3}=1-\zeta_{0}^{2}, \kappa_{2}=\zeta_{3}-\zeta_{1} \zeta_{0}, \kappa_{1}=\zeta_{3}-\zeta_{2} \zeta_{0}, \kappa_{0}=\zeta_{1}-\zeta_{3} \zeta_{0}, \chi_{2}=\kappa_{3}^{2}-\kappa_{0}^{2} \\
\chi_{1}=\kappa_{3} \kappa_{2}-\kappa_{1} \kappa_{0}, \chi_{0}=\kappa_{3} \kappa_{1}-\kappa_{2} \kappa_{0}
\end{array}
\end{equation}

If Nash equilibrium satisfies Eq.(27), it is understood as locally stable. Because Eq. (27) is complex and difficult to be solved, we will further analyze the complexity of the system by performing the numerical simulation.

\section{Complexity analysis}

In this part, the numerical simulation can be utilized to measure the complex dynamic characteristics of dual-channel supply chain. The value of the system parameter is set to be: $k_1=0.1$, $k_2=0.1$, $w_1=0.31$, $w_2=0.4$, $g_1=0.4$, $g_2=0.13$, $g_3=0.12$, $g_4=0.15$, $\mu=0.5$, $\eta=1.1$, $\lambda=0.1$, $g=0.5$, $g_0=0.1$ and $l=0.1$. The stability, bifurcation and chaos behaviors are studied respectively in the following parts.

%Í¼Æ¬1 abcdef
\begin{figure}[!htbp]
\centering

\subfigure[ Model NG, $g_1\in(0,3)$.]{
\begin{minipage}[t]{0.5\linewidth}
\centering
\includegraphics[width=9cm]{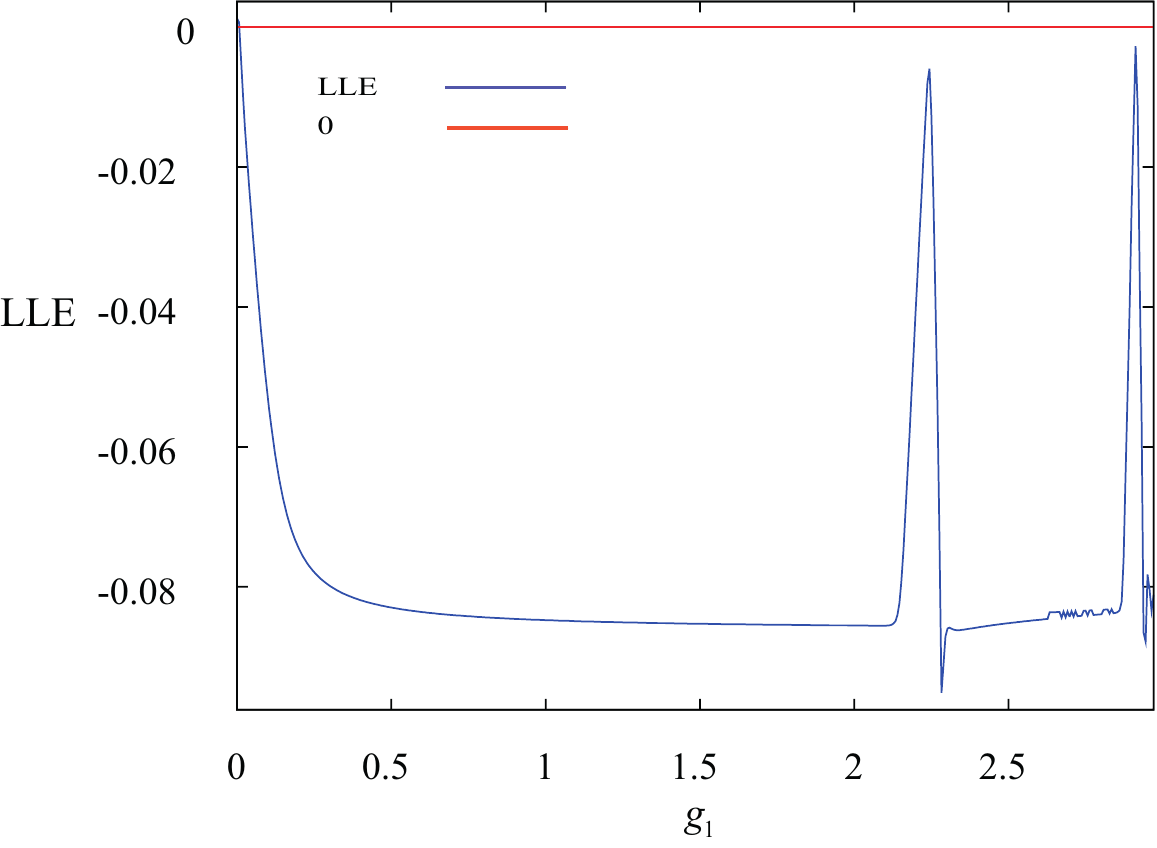}
%\caption{fig1}
\end{minipage}%
}%
\subfigure[ Model NG, $g_1\in(0,3)$.]{
\begin{minipage}[t]{0.5\linewidth}
\centering
\includegraphics[width=9cm]{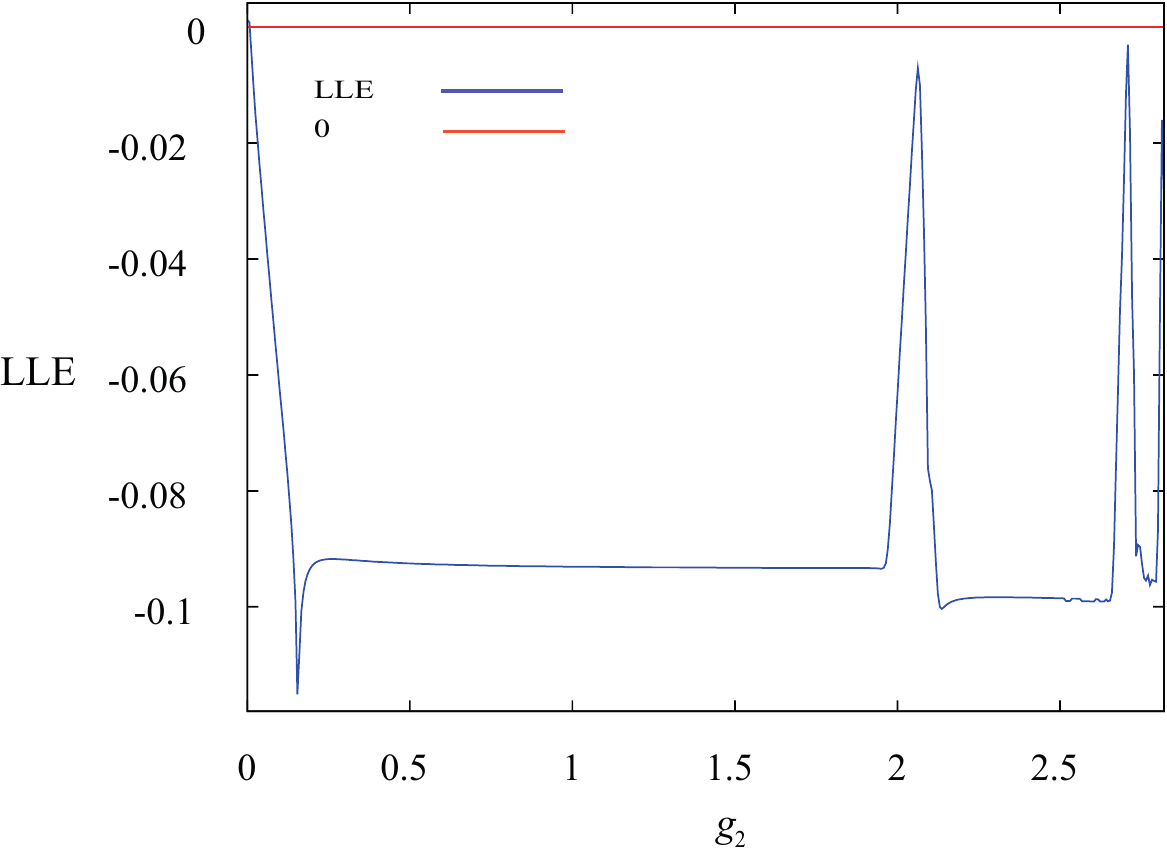}
%\caption{fig2}
\end{minipage}%
}%

\subfigure[ Model MS, $g_1\in(0,3)$.]{
\begin{minipage}[c]{0.5\linewidth}
\centering
\includegraphics[width=9cm]{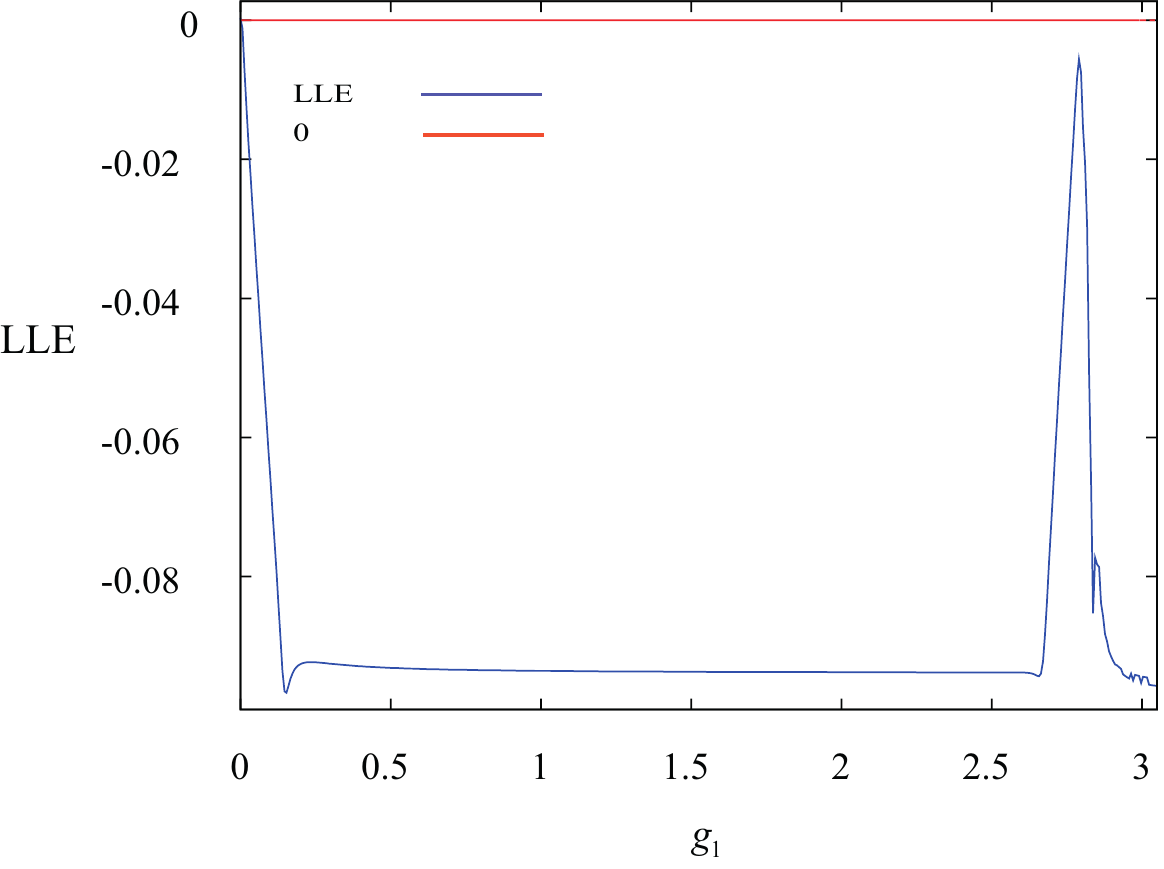}
%\caption{fig2}
\end{minipage}%
}%
\subfigure[ Model MS, $g_1\in(0,3)$.]{
\begin{minipage}[c]{0.5\linewidth}
\centering
\includegraphics[width=9cm]{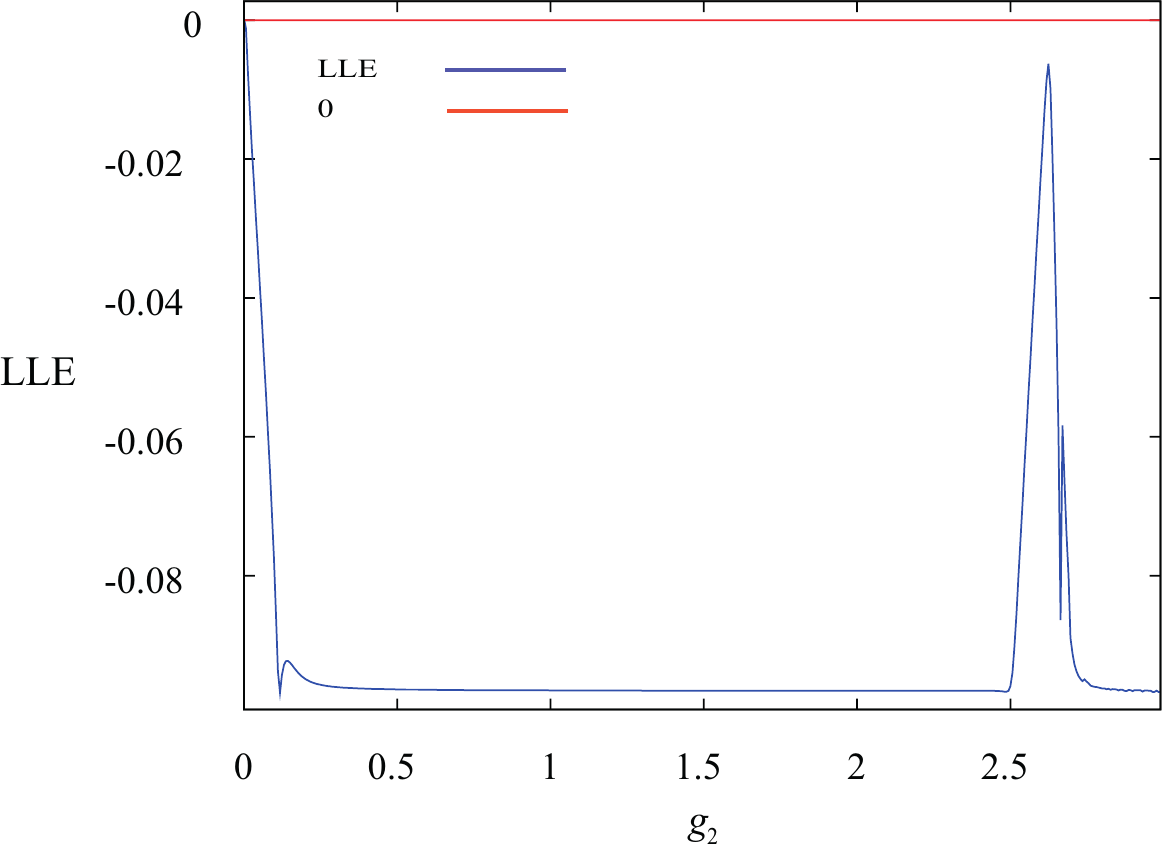}
%\caption{fig2}
\end{minipage}%
}%

\subfigure[ Model RS, $g_1\in(0,4)$.]{
\begin{minipage}[c]{0.5\linewidth}
\centering
\includegraphics[width=9cm]{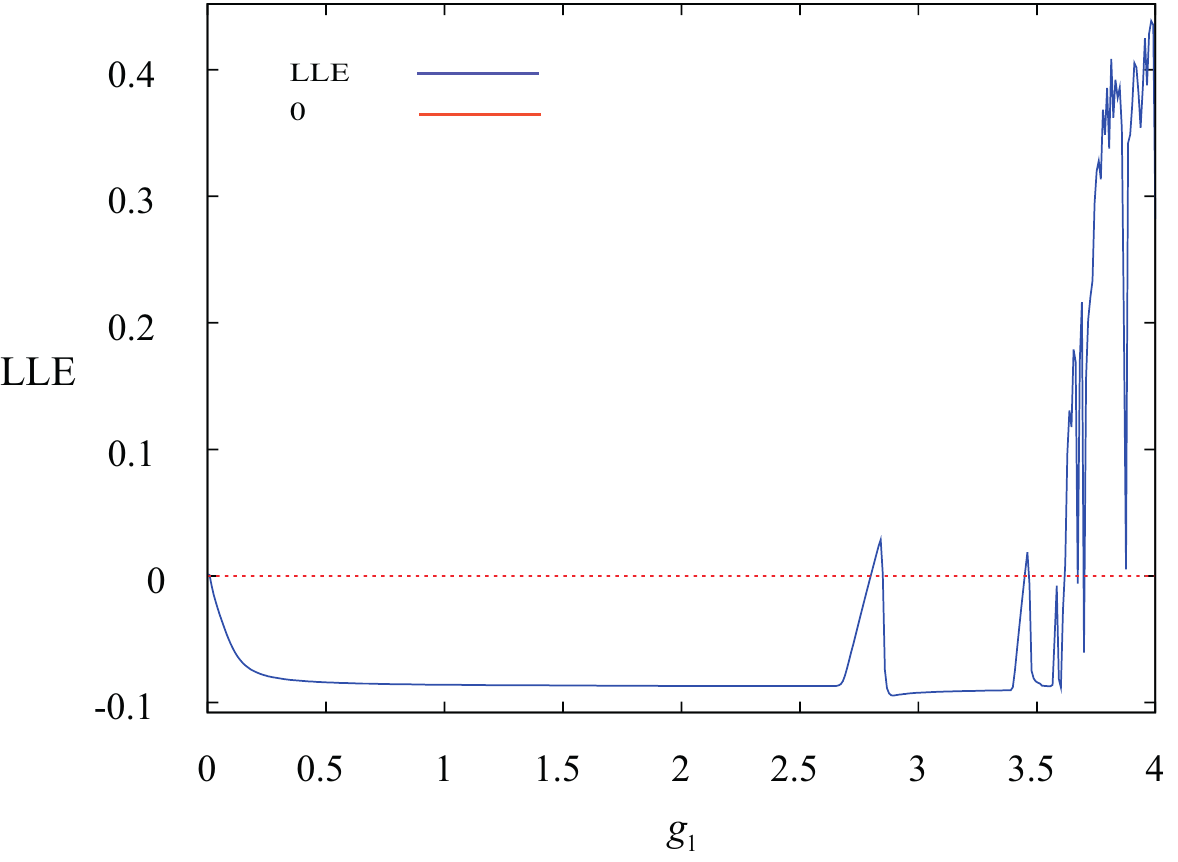}
%\caption{fig2}
\end{minipage}%
}%
\subfigure[ Model RS, $g_1\in(0,4)$.]{
\begin{minipage}[c]{0.5\linewidth}
\centering
\includegraphics[width=9cm]{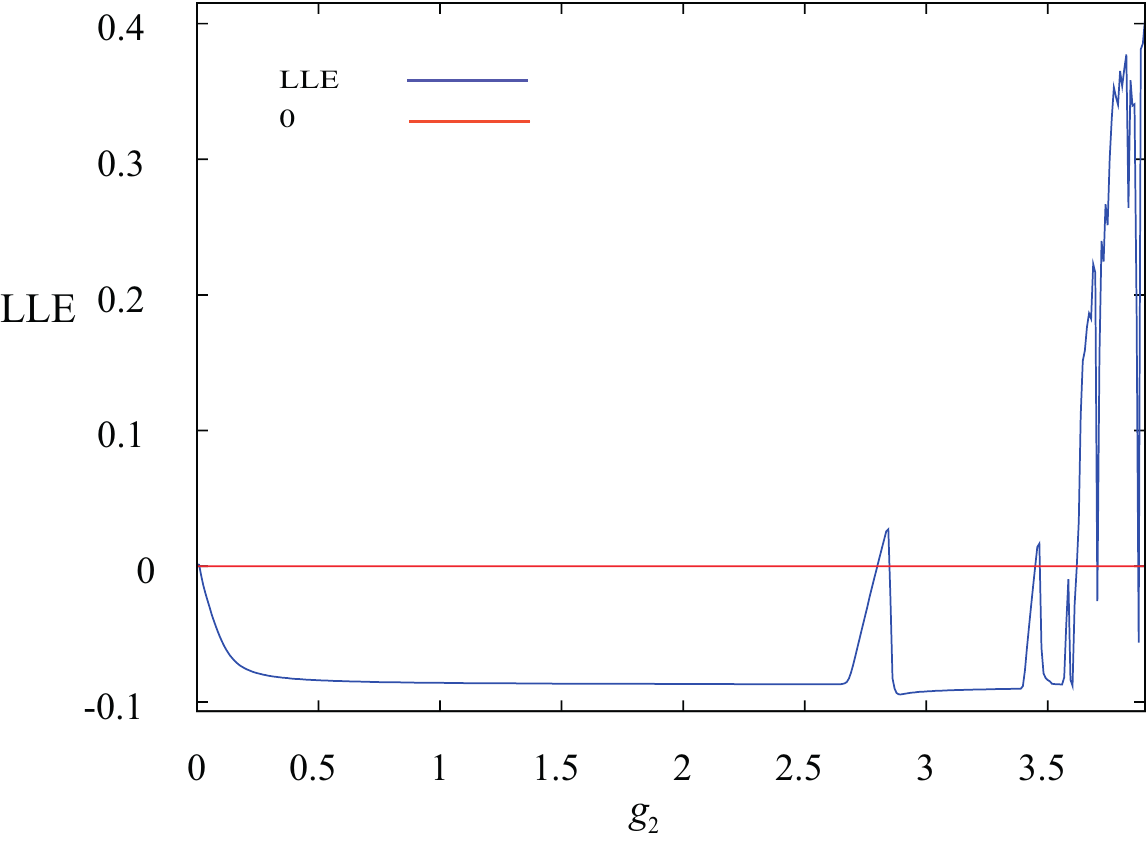}
%\caption{fig2}
\end{minipage}%
}%
\centering
\caption{ The largest Lyapunov exponent with respect to $g_1$ and $g_2$.}
\end{figure}

Figure 1 is the largest Lyapunov exponent diagram with $g_1$,$g_2$ changes under the three models. We can summarize from Figures 1(a), (b) that, under the NG model, the system enters into the four-fold bifurcation region with varied adjustment speeds of sales commissions. Under the MS model, the system only enters into the double bifurcation region with various adjustment speeds of sales commissions in Figures 1(c), (d). Under the RS model, the system finally enters chaos with various adjustment speeds of sales commissions in Figures 1(e), (f). This phenomenon is also verified in the attractor graph below.

%Í¼Æ¬2 abc
\begin{figure}[!htbp]
\centering

\subfigure[ Model NG.]{
\begin{minipage}[t]{0.5\linewidth}
\centering
\includegraphics[width=9cm]{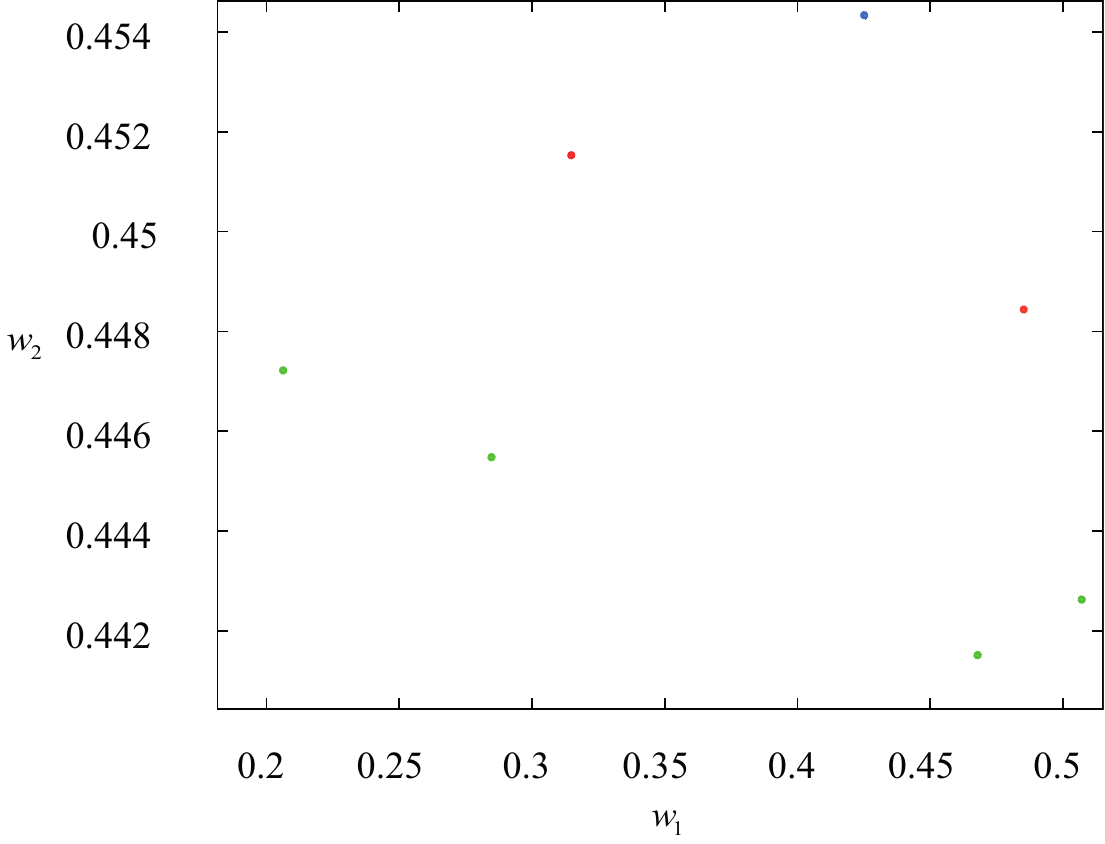}
%\caption{fig1}
\end{minipage}%
}%
\subfigure[ Model MS.]{
\begin{minipage}[t]{0.5\linewidth}
\centering
\includegraphics[width=9cm]{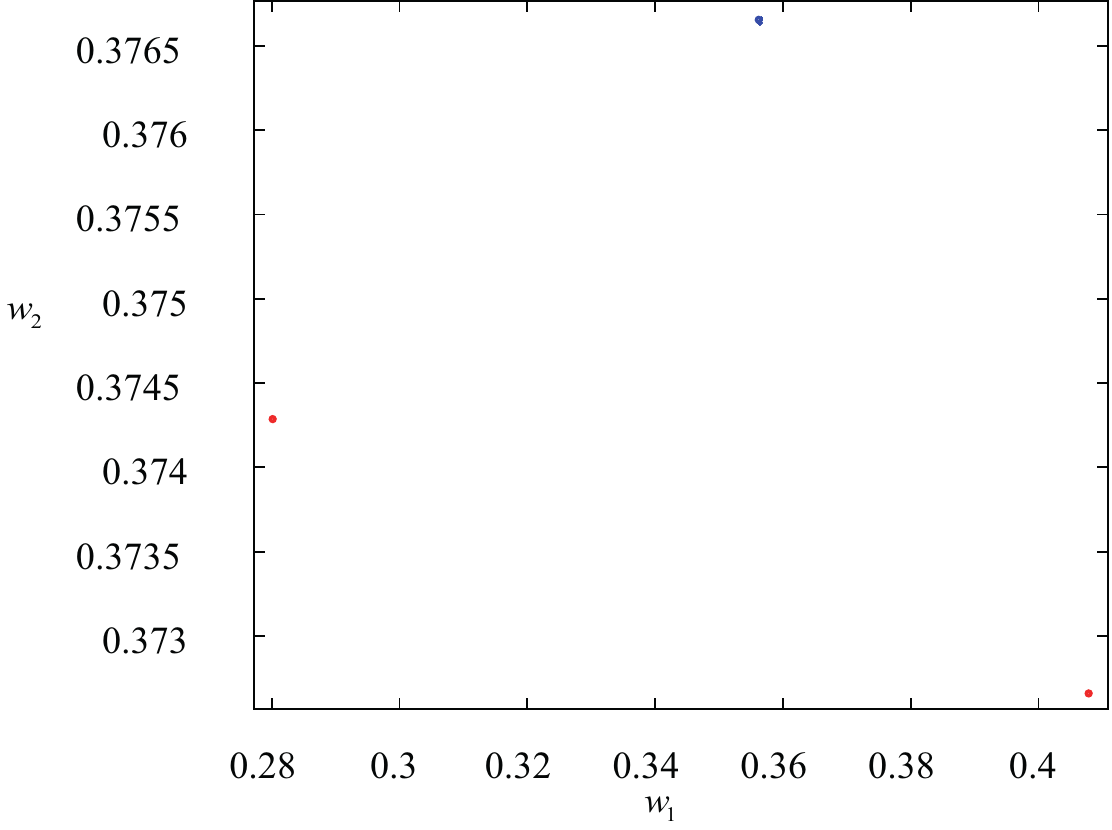}
%\caption{fig2}
\end{minipage}%
}%

\subfigure[ Model RS.]{
\begin{minipage}[c]{0.5\linewidth}
\centering
\includegraphics[width=9cm]{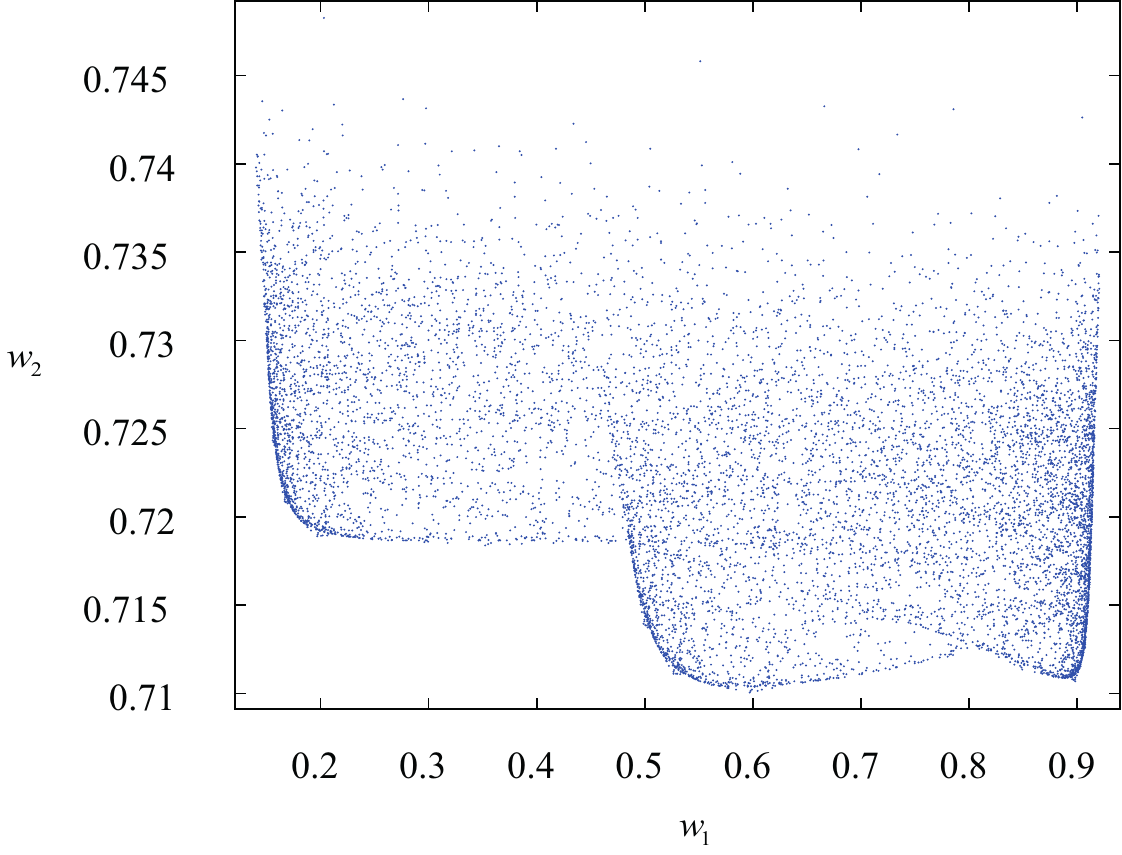}
%\caption{fig2}
\end{minipage}%
}%
\centering
\caption{The chaotic attractor under three power structures.}
\end{figure}

Figure 2 is plotted to explain the attractor graphs of the NG, MS and RS models. In Figure 2(a), the model NG finally enters the 4-period point. It can be seen that the system gradually transits from a blue fixed point to two red 2-period points and four green 4-period points. In Figure 2(b), the system only enters into the double bifurcation region, so only one blue fixed point and two red 2-period points appear in model MS. Figure 2(c) shows that in the model RS, when $g_1=4$, the system finally enters into the chaotic state.

%Í¼Æ¬3 abcdef

\begin{figure}[!htbp]
\centering

\subfigure[ Model NG, $g_1=0.4$.]{
\begin{minipage}[t]{0.5\linewidth}
\centering
\includegraphics[width=7cm]{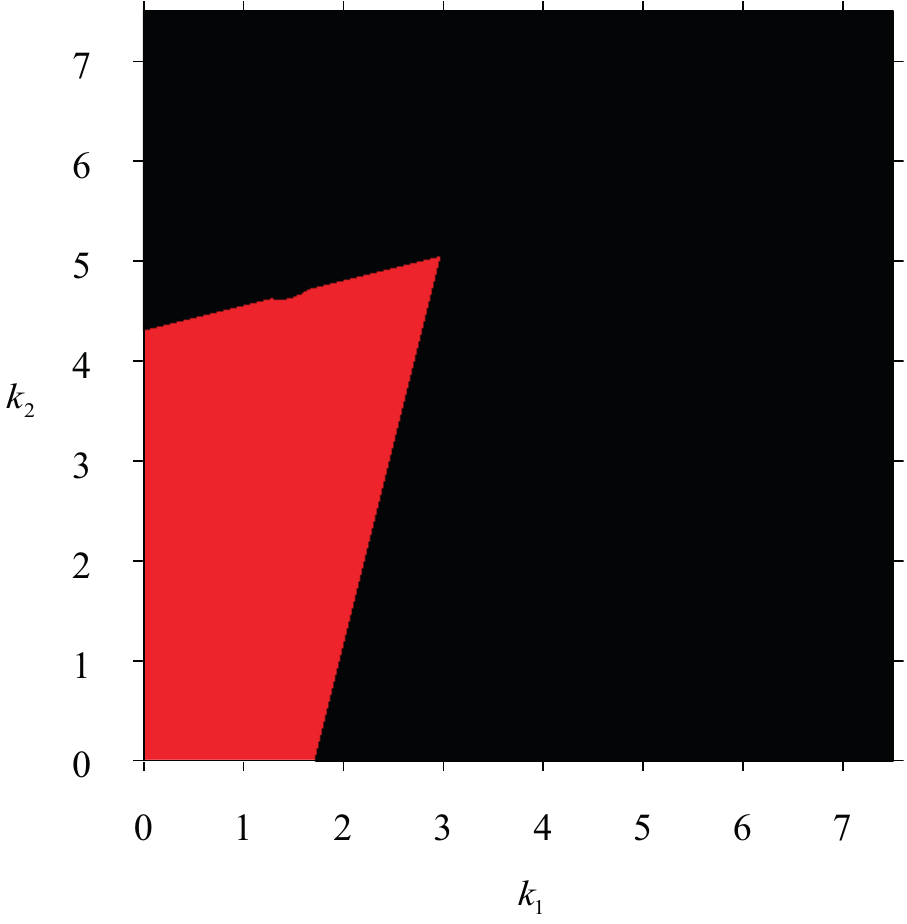}
%\caption{fig1}
\end{minipage}%
}%
\subfigure[ Model NG, $g_1=2.1$.]{
\begin{minipage}[t]{0.5\linewidth}
\centering
\includegraphics[width=7cm]{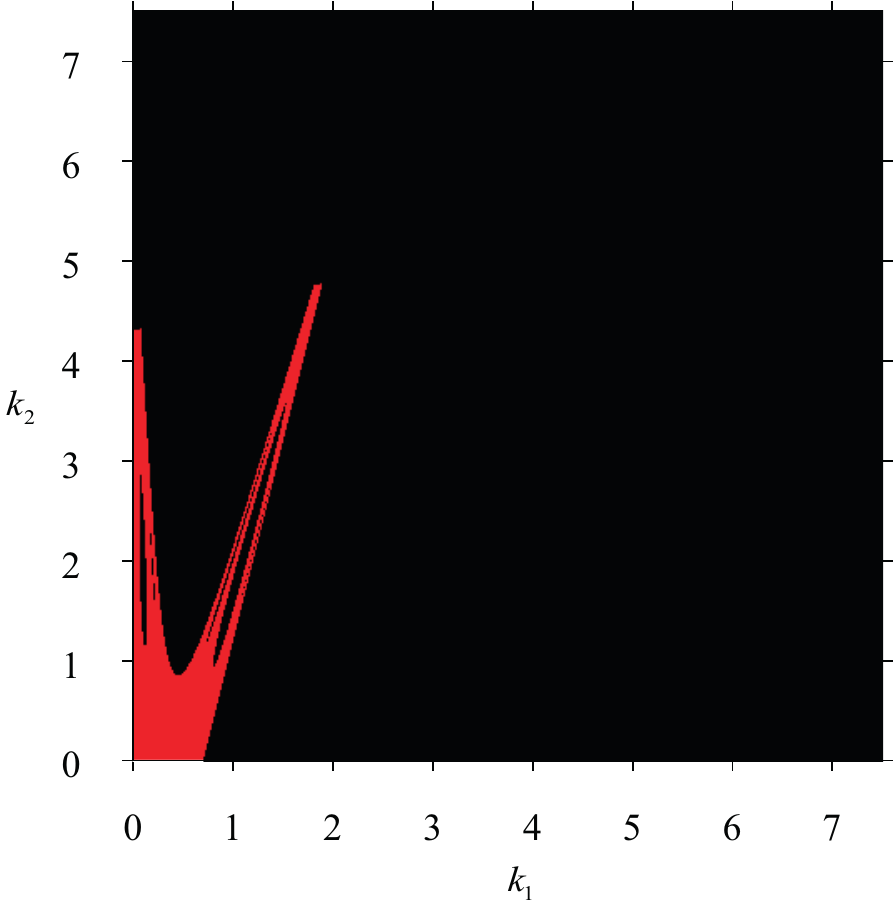}
%\caption{fig2}
\end{minipage}%
}%

\subfigure[ Model MS, $g_1=0.4$.]{
\begin{minipage}[c]{0.5\linewidth}
\centering
\includegraphics[width=7cm]{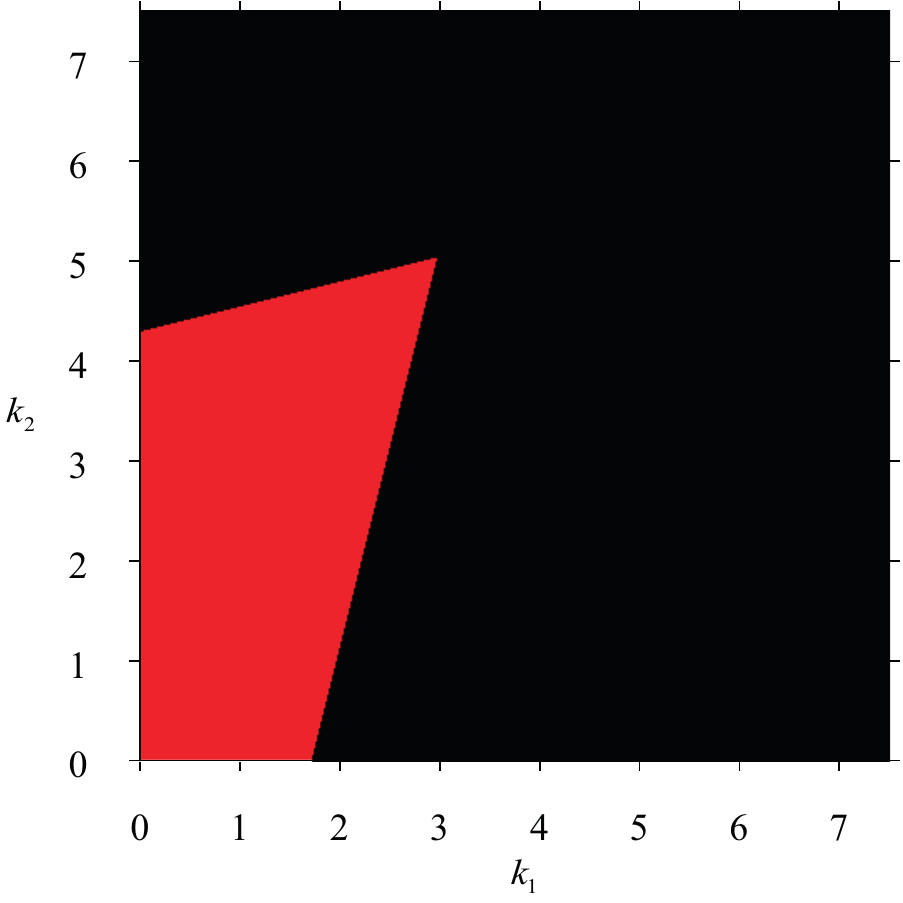}
%\caption{fig2}
\end{minipage}%
}%
\subfigure[ Model MS, $g_1=2$.]{
\begin{minipage}[c]{0.5\linewidth}
\centering
\includegraphics[width=7cm]{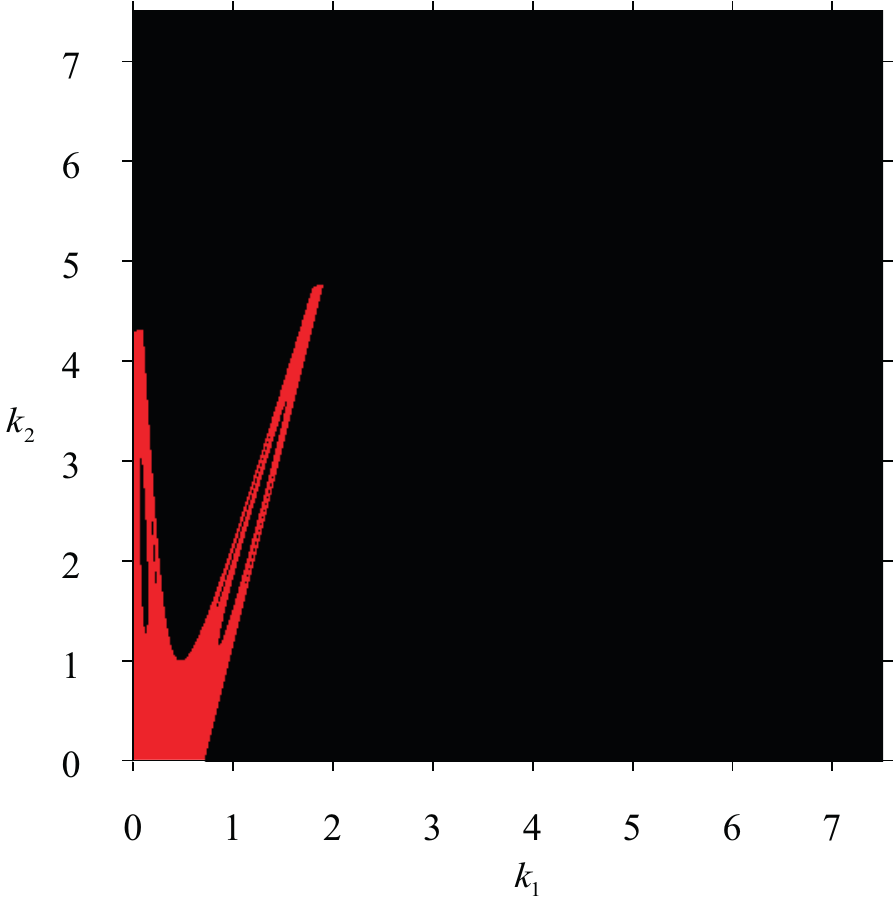}
%\caption{fig2}
\end{minipage}%
}%

\subfigure[ Model RS, $g_1=0.4$.]{
\begin{minipage}[c]{0.5\linewidth}
\centering
\includegraphics[width=7cm]{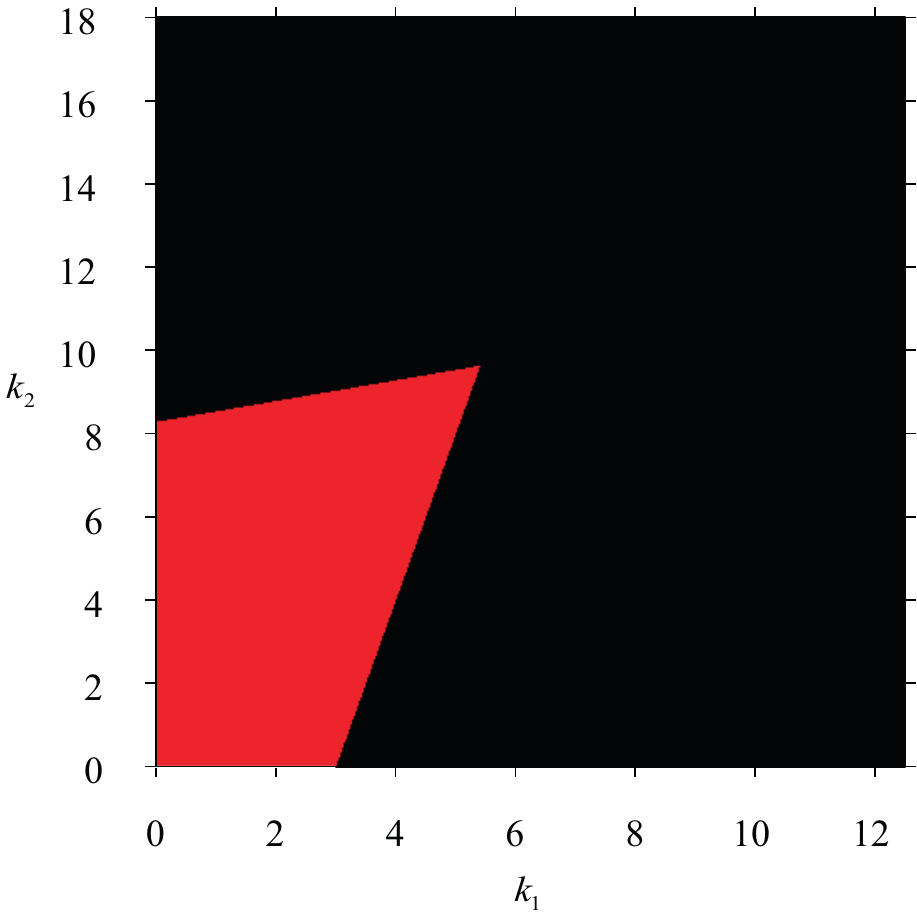}
%\caption{fig2}
\end{minipage}%
}%
\subfigure[ Model RS, $g_1=2.7$.]{
\begin{minipage}[c]{0.5\linewidth}
\centering
\includegraphics[width=7cm]{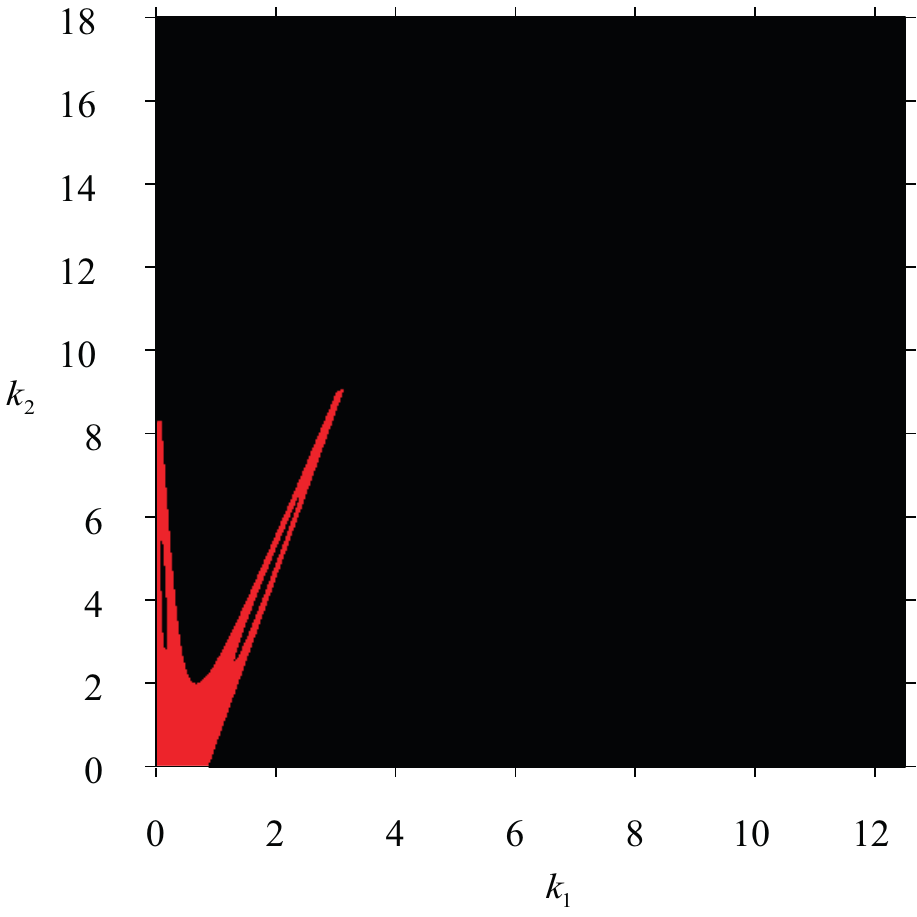}
%\caption{fig2}
\end{minipage}%
}%
\centering
\caption{Attraction domain in stable and unstable regions under three power structures.}
\end{figure}

Figure 3 shows the variation of the attraction domain under the three models. In Figures 3(a), 3(c) and 3(e), the red area indicates the attraction domain for the stable state. Figures 3(b), 3(d) and 3(f) show attraction domains when $g_1=2.1$, $g_1=2$ and $g_1=2.7$ in the model NG, MS, and RS respectively. We found that when the system is in an unstable area, there are some missing parts in the red region, which means online and offline retailers have reduced their regions of attraction of sales commission.

%Í¼Æ¬4 abc

\begin{figure}[!htbp]
\centering

\subfigure[ Model NG.]{
\begin{minipage}[t]{0.5\linewidth}
\centering
\includegraphics[width=9cm]{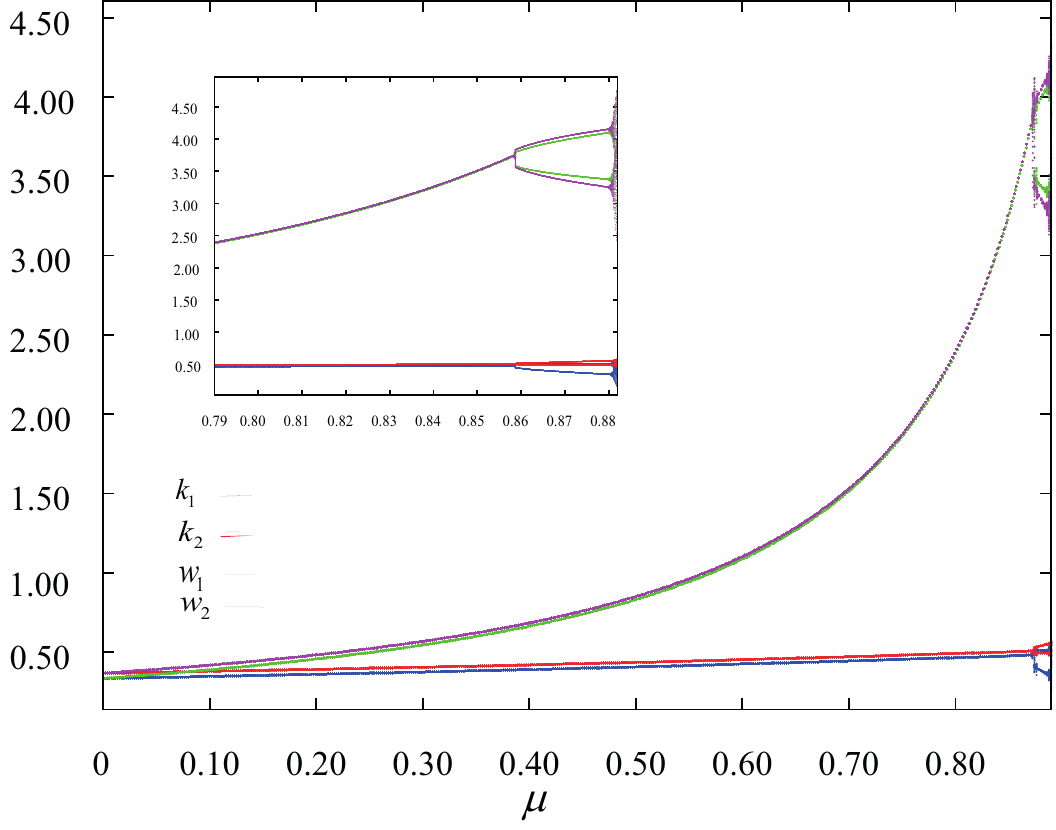}
%\caption{fig1}
\end{minipage}%
}%
\subfigure[ Model MS.]{
\begin{minipage}[t]{0.5\linewidth}
\centering
\includegraphics[width=9cm]{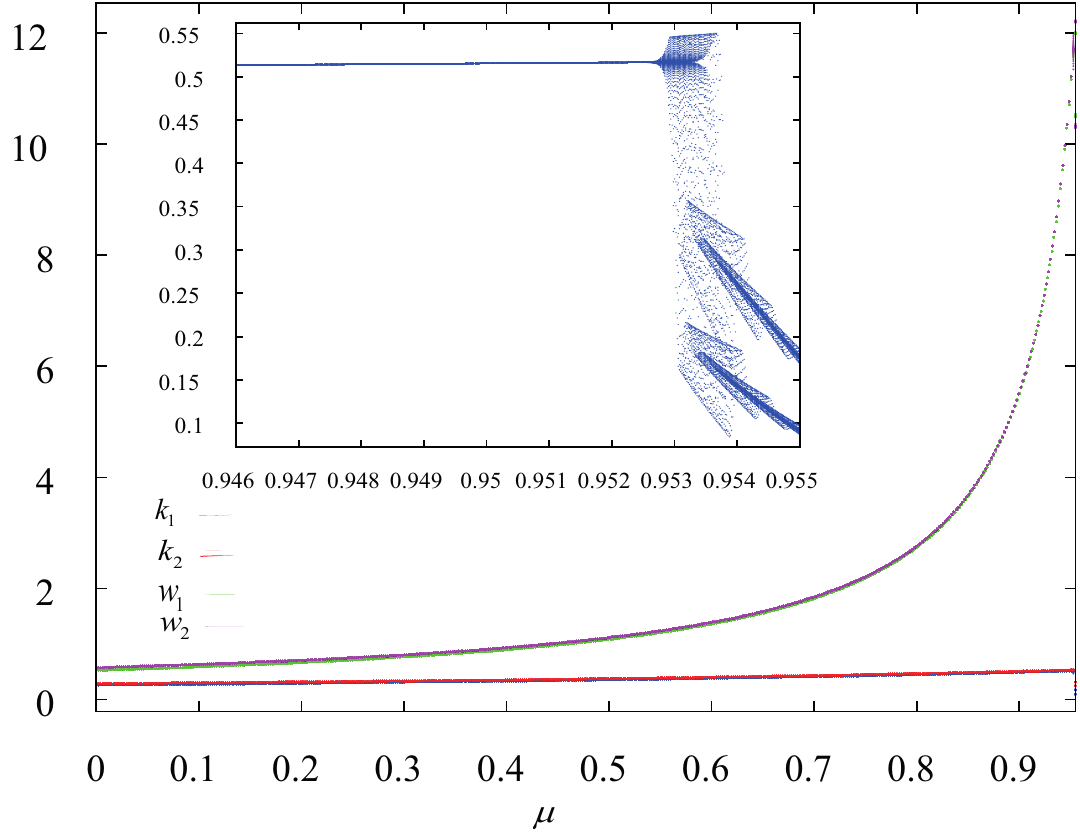}
%\caption{fig2}
\end{minipage}%
}%

\subfigure[ Model RS.]{
\begin{minipage}[c]{0.5\linewidth}
\centering
\includegraphics[width=9cm]{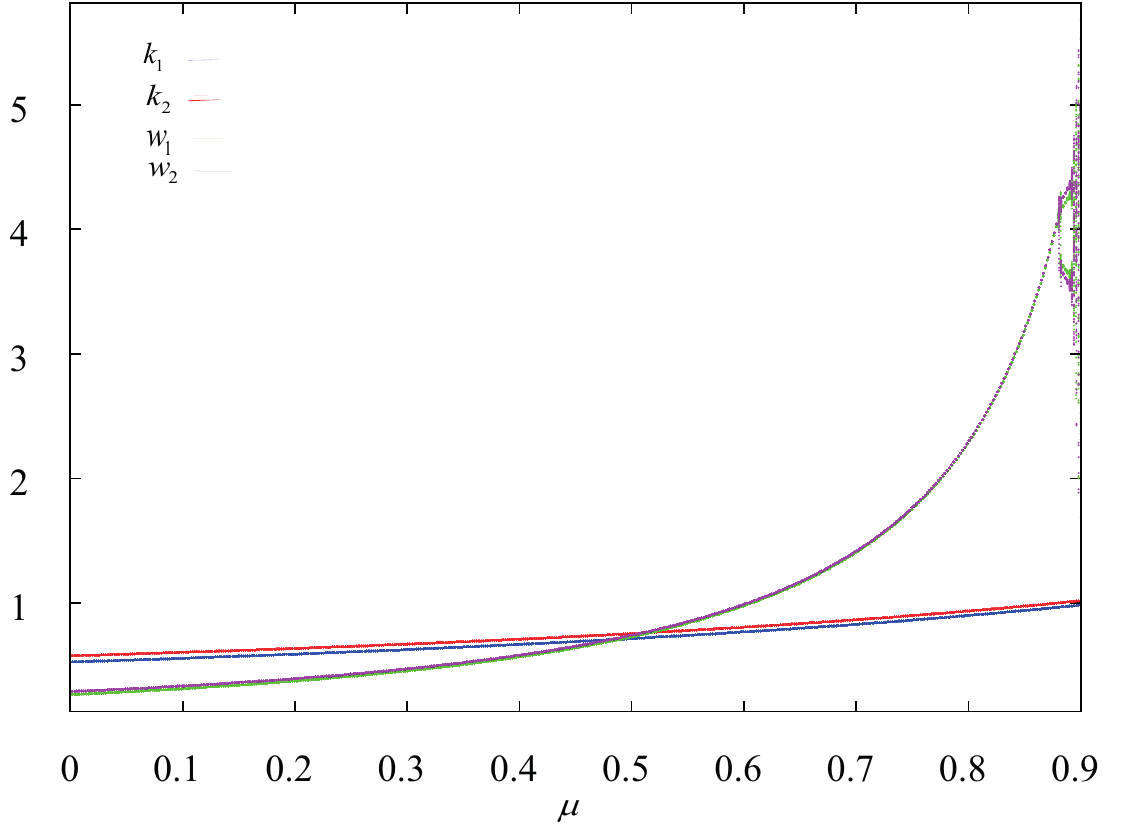}
%\caption{fig2}
\end{minipage}%
}%
\centering
\caption{The bifurcation diagram under different scenarios w.r.t $\mu$.}
\end{figure}

Figure 4 respectively depicts the bifurcation diagram with respect to $\mu$ under the three models. It shows the complexity of the retailer's substitution. In all bifurcation diagrams in this paper, the blue line stands for $k_1$ and the red line stands for $k_2$. The green line represents $w_1$, and the purple line represents $w_2$. $\mu=0$ represents the complete independence of demand, and retailers are monopolists in their regions. And  $\mu=1$ represents that the two retailers are completely replaceable.

Figure 4(a) shows the changing tendency of $\mu$in the model NG. $k_1$ and $k_2$, $w_1$ and $w_2$ all gradually increase, but the growth of $k_1$ and $k_2$ never exceeds that of $w_1$ and $w_2$ . At $\mu=0.86$, the price system begins the 2-period bifurcation. When $\mu$ keeps on increasing change by degrees, the price system suffers greater volatility. We can conclude that the pricing system of the manufacturer and retailers enters an unstable area due to the increase in retailer substitutability. In the small picture located in the upper left corner of Figure 4(b), it can be clearly seen that in the model MS, when $\mu=0.952$, the system directly enters chaos. This is later than the first bifurcation point of the model NG. Therefore, it's evident that in the model MS, the retailer's substitution effect is weaker and the system is more stable. According to Figure 4(c), in the model RS, when $\mu=0.49$, $w_1$ and $w_2$ exceed $k_1$ and $k_2$. When $\mu=0.86$, bifurcation of $w_1$ and $w_2$ occurs, and $k_1$ and $k_2$ always remain stable. Retailers dominate the supply chain system, making the two retailers' sales commissions more profitable at the beginning. However, with the fierce market competition and the increase in the two retailers' substitutability, the manufacturer's wholesale prices eventually exceed the retailers' sales commissions. At the same time, retailers still rely on their dominant position to keep their sales commissions stable, but manufacturer's wholesale prices will fall into the chaotic region due to fierce market competition.

In summary, no matter what kind of power structure model is used, the manufacturer will get more attention to the substitutability of retailers. And the increase in the substitutability of retailers leads to the rise in wholesale prices. However, substitutability between retailers has little effect on retailers' sales commission. As a result, fierce competition between retailers may lead manufacturer to raise wholesale prices. As long as there is a certain degree of substitutability between retailers of different channels, the more intense the competition, the more instability will happen in the real market. It is worth noting that under the manufacturer-dominated supply chain model, the price system enters into chaos relatively late, so it is more stable.

%Í¼Æ¬5 abc

\begin{figure}[!htbp]
\centering

\subfigure[ Model NG.]{
\begin{minipage}[t]{0.5\linewidth}
\centering
\includegraphics[width=9cm]{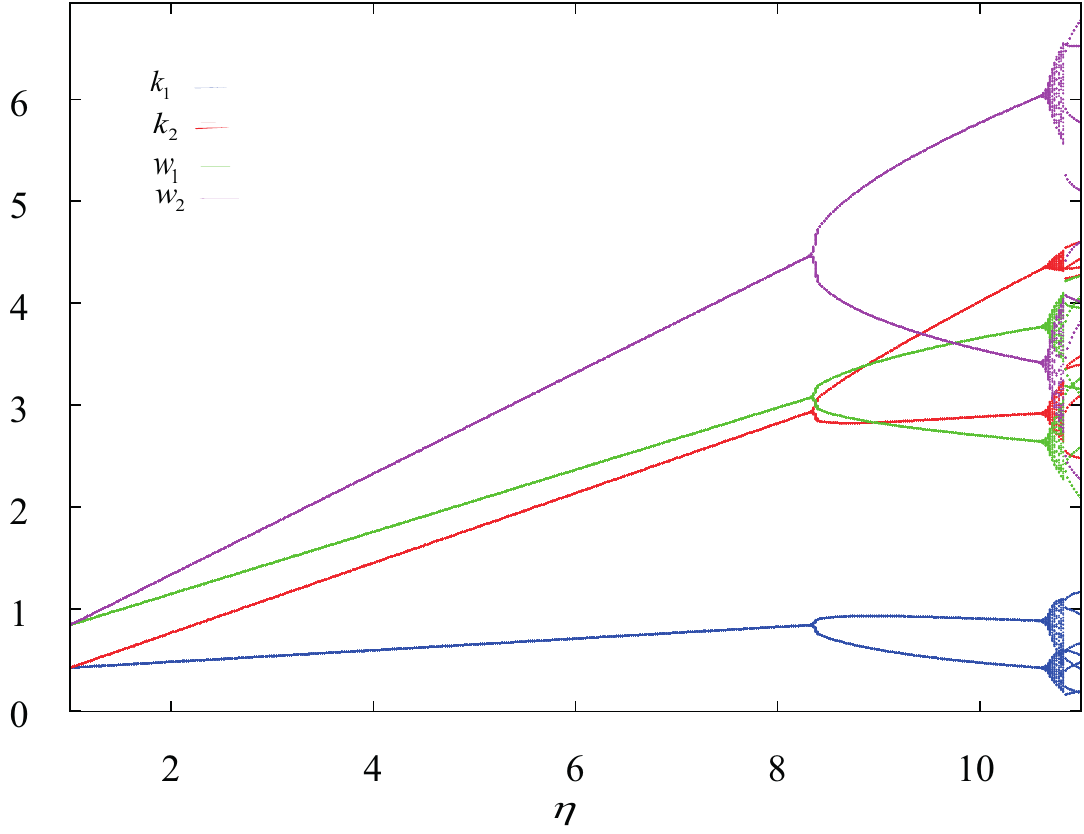}
%\caption{fig1}
\end{minipage}%
}%
\subfigure[ Model MS.]{
\begin{minipage}[t]{0.5\linewidth}
\centering
\includegraphics[width=9cm]{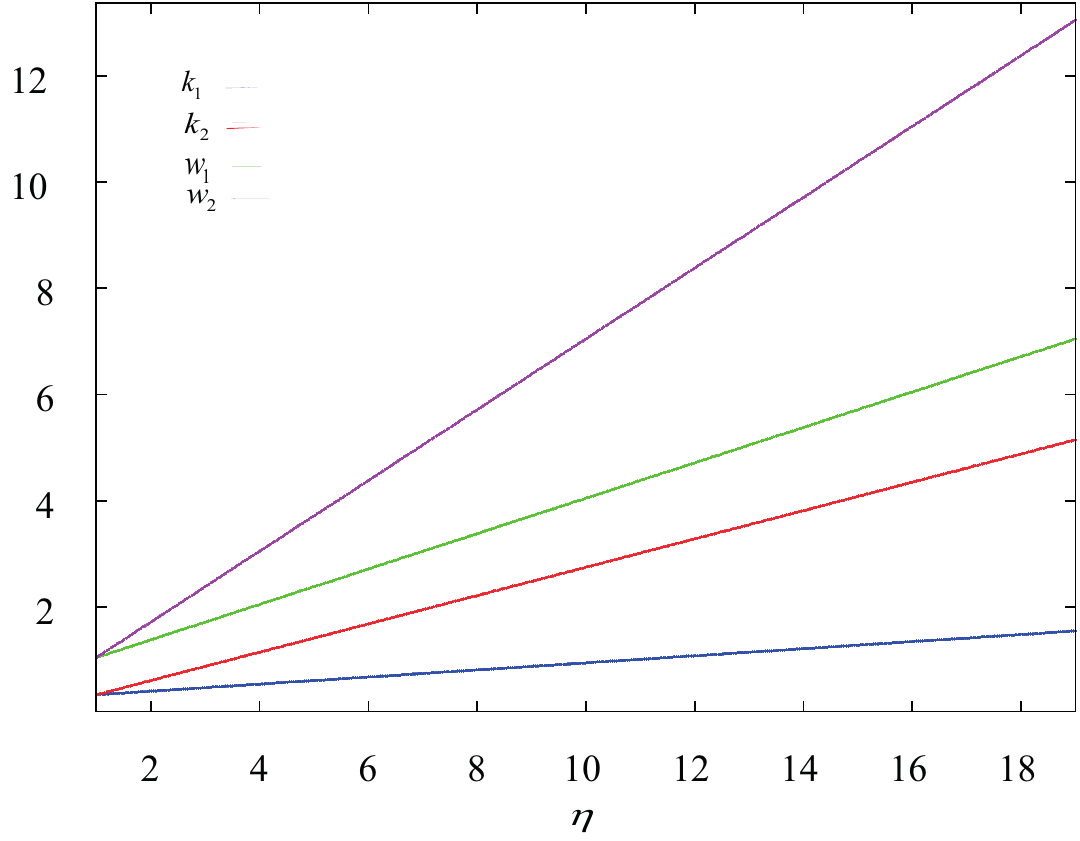}
%\caption{fig2}
\end{minipage}%
}%

\subfigure[ Model RS.]{
\begin{minipage}[c]{0.5\linewidth}
\centering
\includegraphics[width=9cm]{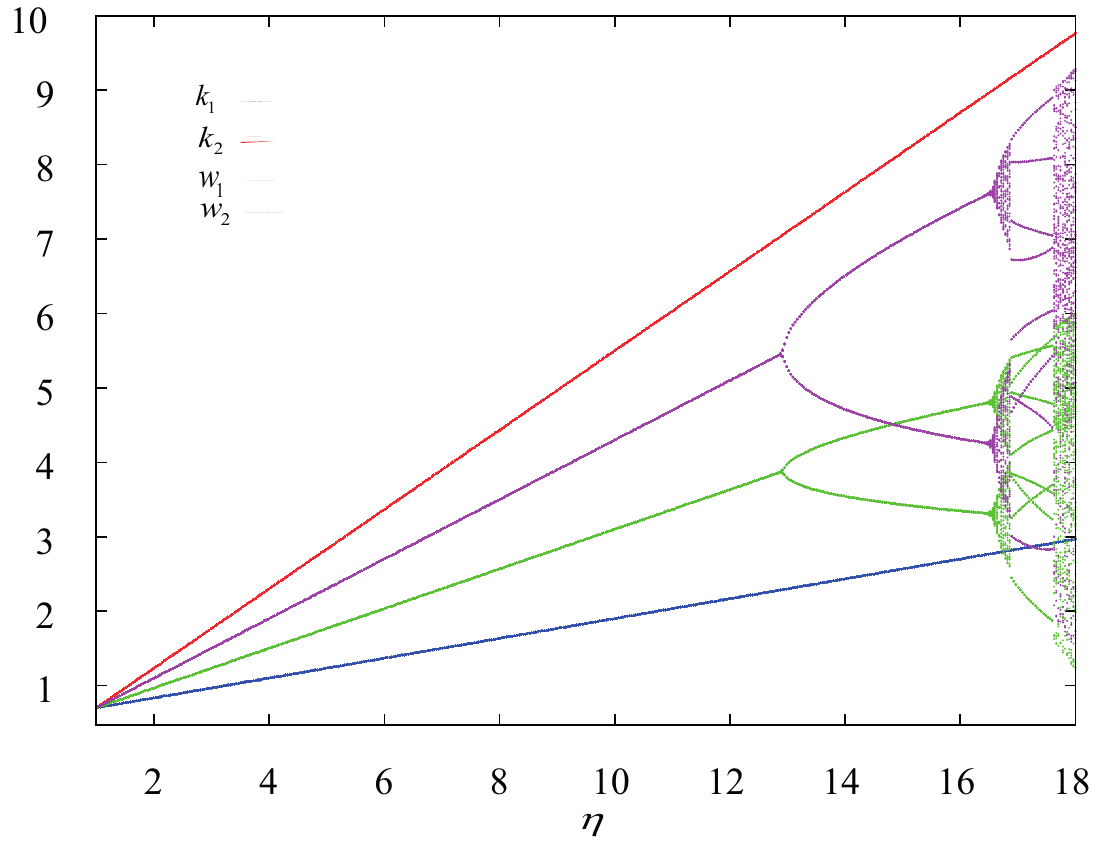}
%\caption{fig2}
\end{minipage}%
}%
\centering
\caption{The bifurcation diagram under different scenarios w.r.t $\eta$.}
\end{figure}

Next, Figure 5 describes the complex characteristics as consumers' channel preference $\eta$ varies. Figure 5 (a) depicts that in the model NG, as $\eta$ increasing, the system gradually falls into the chaotic region. As shown in the figure, if the online channel keeps a greater degree of consumers' preferences, wholesale price and sales commission will increase, where $ w_2>w_1$ and $k_2>k_1$. However, if consumers excessively prefer online channel, the system may appear bifurcation and chaos. Figure 5 (b) indicates that in the model MS, with the growth of consumers' channel preference, the price system does not appear bifurcation, and it is relatively stable. As seen in Figure 5 (c), there exists that in the model RS, when $\eta$ increases to 12.96, the manufacturer's wholesale price bifurcates and then gradually falls into chaos, but the retailers' sales commissions do not bifurcate. Retailers dominate the market and consumers' channel preference is too large. It will lead to dramatic fluctuations in manufacturer's wholesale prices.

%Í¼Æ¬6 abc
\begin{figure}[!htbp]
\centering

\subfigure[ Model NG.]{
\begin{minipage}[t]{0.5\linewidth}
\centering
\includegraphics[width=9cm]{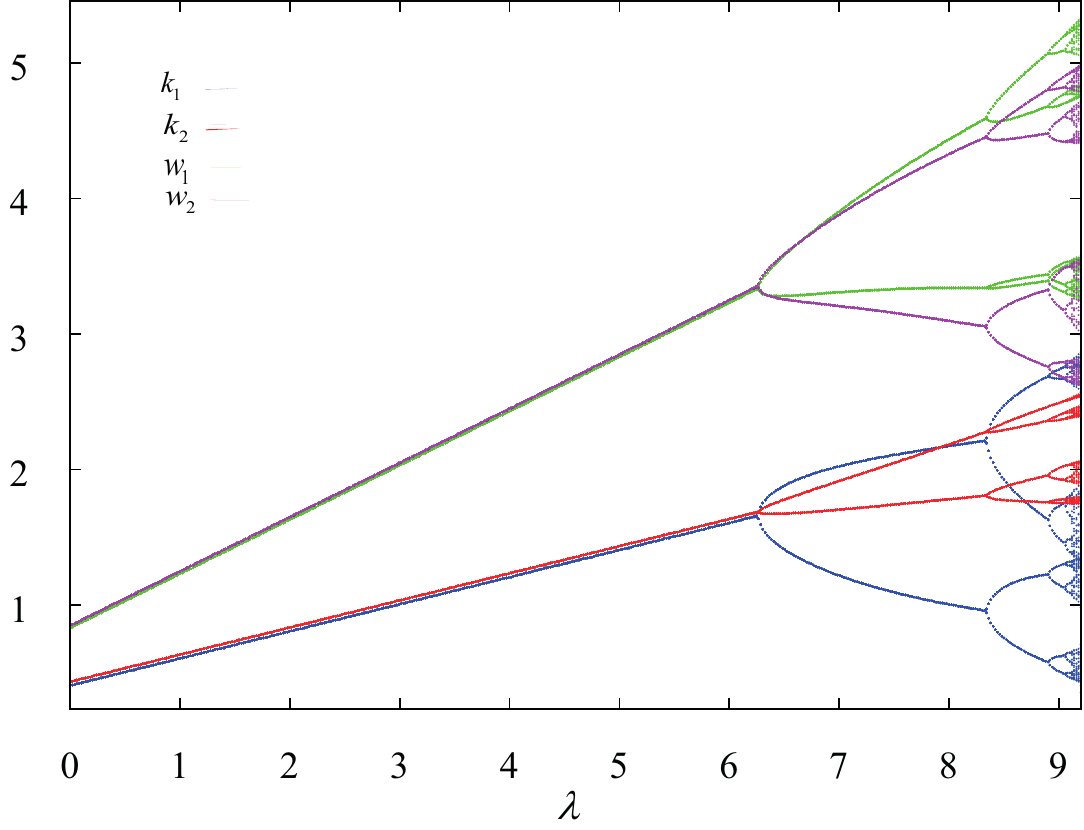}
%\caption{fig1}
\end{minipage}%
}%
\subfigure[ Model MS.]{
\begin{minipage}[t]{0.5\linewidth}
\centering
\includegraphics[width=9cm]{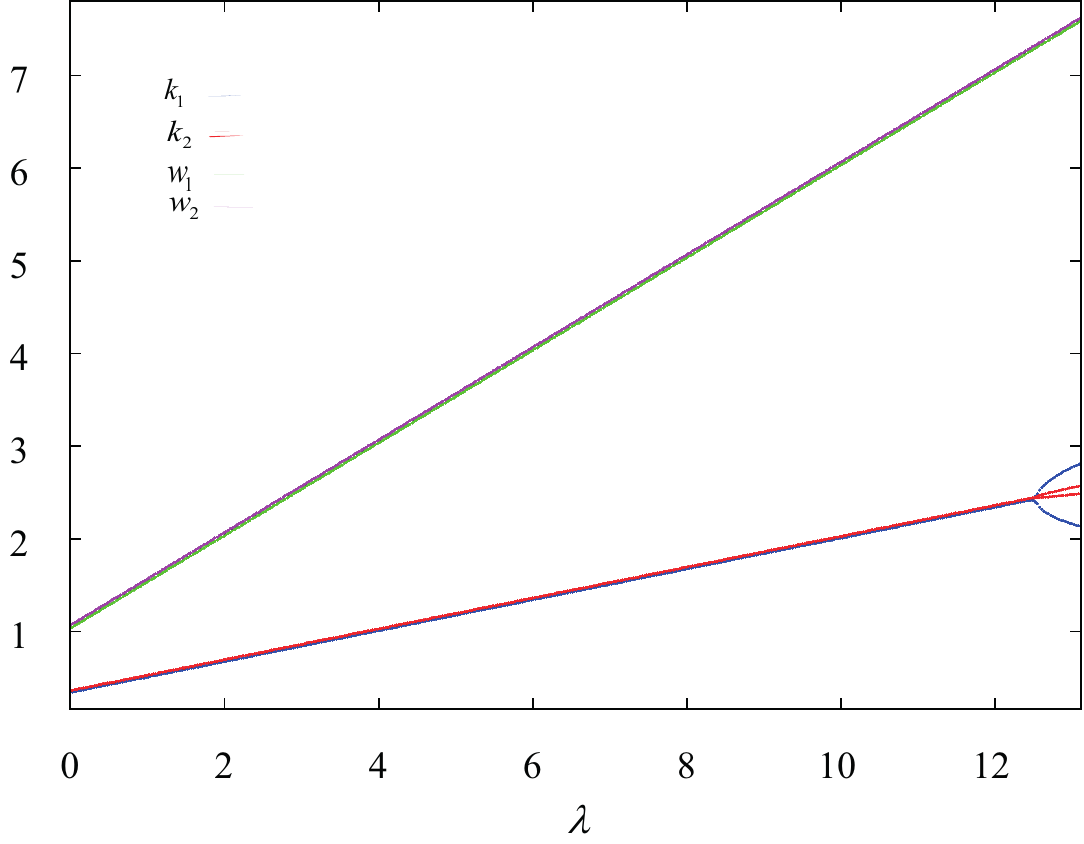}
%\caption{fig2}
\end{minipage}%
}%

\subfigure[ Model RS.]{
\begin{minipage}[c]{0.5\linewidth}
\centering
\includegraphics[width=9cm]{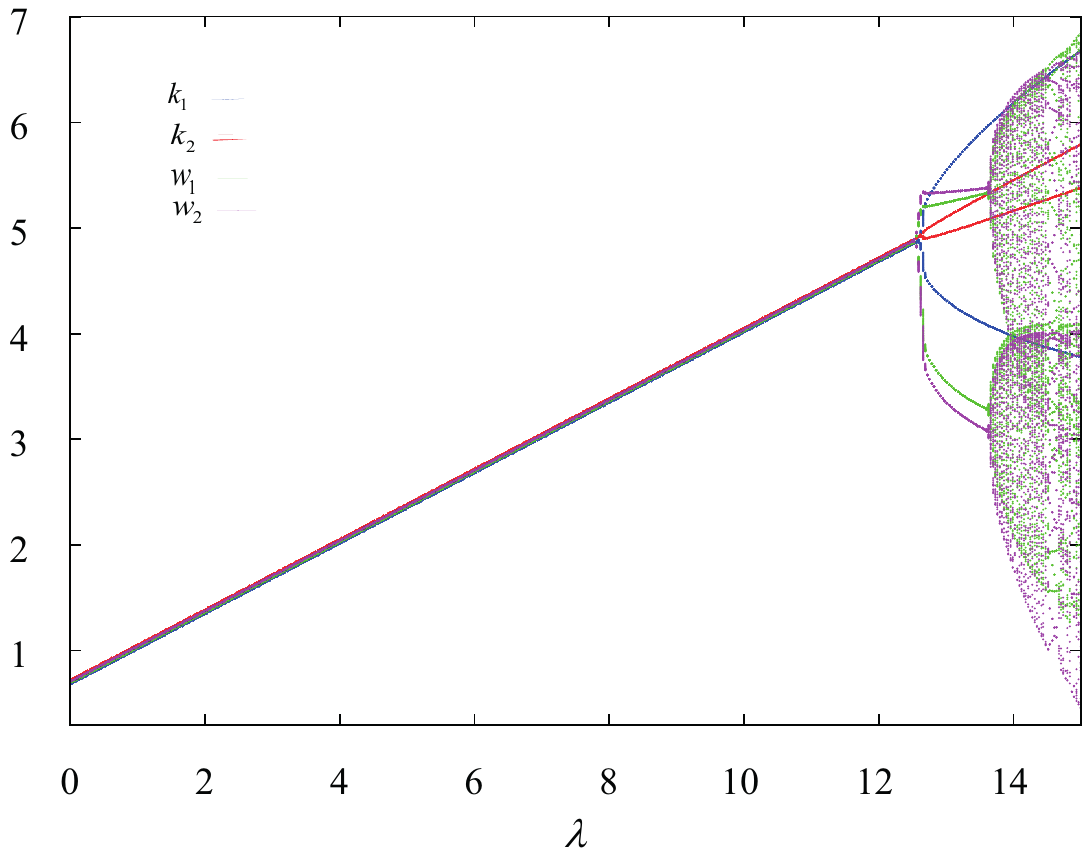}
%\caption{fig2}
\end{minipage}%
}%
\centering
\caption{The bifurcation diagram under different scenarios w.r.t $\lambda$.}
\end{figure}
Figure 6 depicts the bifurcation diagram of consumers' low-carbon preference under three models. In model NG, with the growth of $\lambda$, the manufacturer's wholesale prices $w_1,w_2$ and the two retailers' sales commissions $k_1,k_2$ have increased significantly, as shown in Figure 6(a). Besides, $w_1$ almost coincides with $w_2$, and it is also the same case for $k_1$ and $k_2$. It indicates that $\lambda$ has nearly the same impact on the wholesale price and sales commissions of different channels. When $\lambda=6.17, 8.03, 8.79$ and $9.03$, the system enters the 2, 4 and 8-fold cycle and chaotic state, respectively. This indicates that when consumers pay appropriate attention to the low-carbon attributes of products, manufacturer and retailers can take the opportunity to increase wholesale prices and sales commissions. However, consumers too much focus on the green degree, which will harm the stability of the pricing system of manufacturer and retailers. Figure 6 (b) gives that in the model MS, with the growth of $\lambda$, the retailers' sales commissions bifurcated at $\lambda=12.41$. As shown in Figure 6 (b), when the manufacturer dominated the market, consumers' low-carbon preferences may have an adverse impact on retailers' sales commission. As seen in Figure 6 (c), in the model RS, with the growth of low-carbon preferences, the wholesale price and sales commission almost coincide at the beginning. At $\lambda=12.6$, the bifurcation of the system occurs. With the continuous growth of $\lambda$, the price system finally enters into the chaotic region after multiple iterations. This indicates that retailers dominate the market and the increase in consumers' low-carbon preference causes retailers to set higher sales commissions. And they will set sales commissions that are almost the same as the wholesale price. However, excessive consumers' low-carbon preferences eventually also lead the system into period-doubling bifurcation and chaos.

%Í¼Æ¬7 abc

\begin{figure}[!htbp]
\centering

\subfigure[ Model NG.]{
\begin{minipage}[t]{0.5\linewidth}
\centering
\includegraphics[width=9cm]{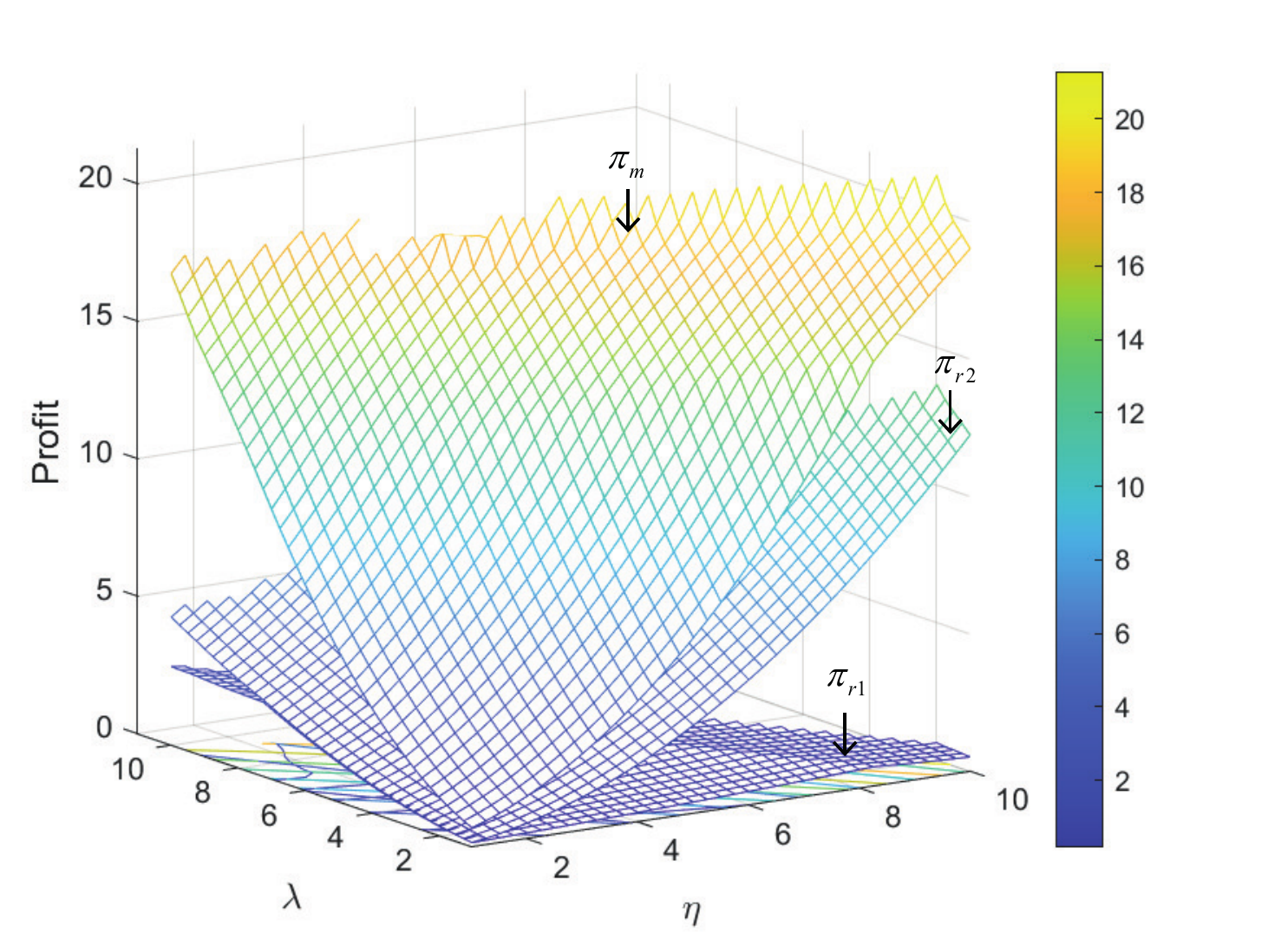}
%\caption{fig1}
\end{minipage}%
}%
\subfigure[ Model MS.]{
\begin{minipage}[t]{0.5\linewidth}
\centering
\includegraphics[width=9cm]{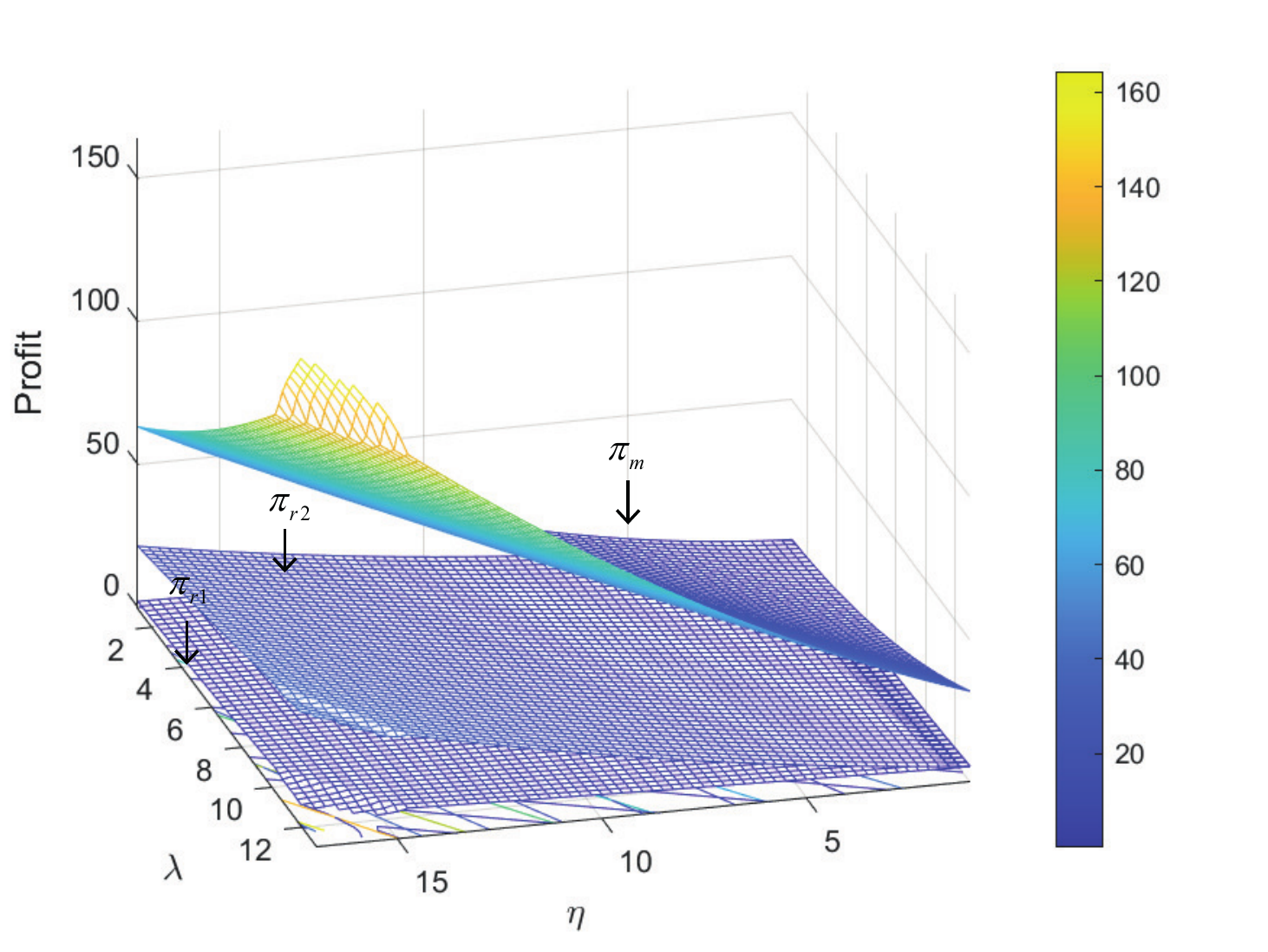}
%\caption{fig2}
\end{minipage}%
}%

\subfigure[ Model RS.]{
\begin{minipage}[c]{0.5\linewidth}
\centering
\includegraphics[width=9cm]{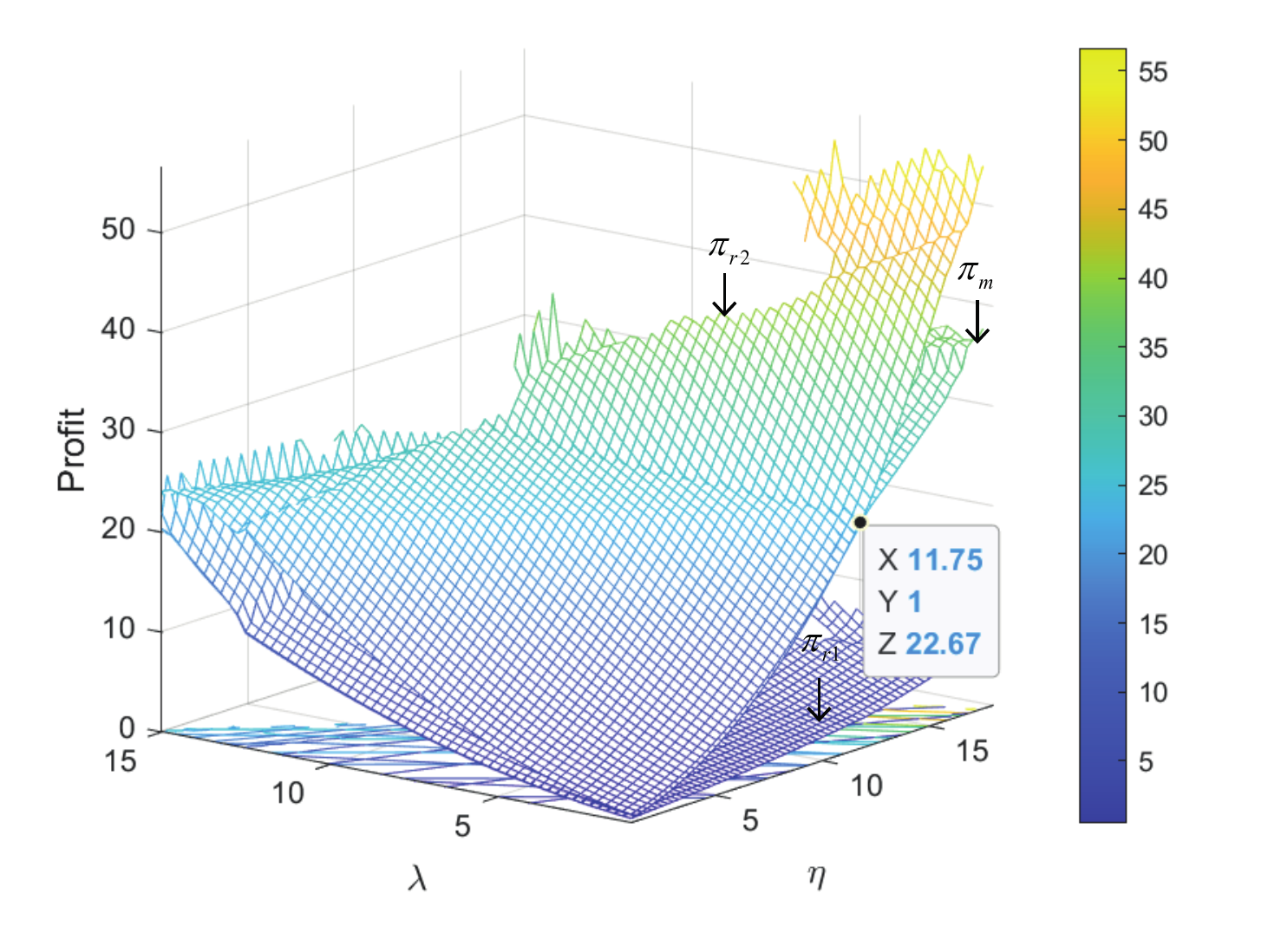}
%\caption{fig2}
\end{minipage}%
}%
\centering
\caption{The average profits of the players under three models.}
\end{figure}

Figure 7 shows the 3D average profits of the manufacturer and two retailers under variations of $\eta$ and $\lambda$ in three models. From Figure 7(a), in model NG, with the growth of channel preference $\eta$, the online channel profit represented by retailer 2 increases more significantly. The increase of  $\lambda$ will moderately boost the profits of the three players. However, at the intersection of $\eta$ and $\lambda$, the average profits of all parties fluctuate sharply and arise missing. Figure 7(b) represents the same profit growth trend of manufacturer and retailers in model MS as Figure 7(a). The difference is that with the growth of low-carbon preference, the system only enters a 2-fold bifurcation, so it can be seen from the color bar on the right that the profits of all parties are greater. And only at the intersection of $\eta$ and $\lambda$, there is only a small fluctuation in the profits of each party. Figure 7(c) indicates that in model RS, with the growth of channel preference, the profits of retailer 2 exceed that of the manufacturer when $\eta=11.75$. At the intersection of $\eta$ and $\lambda$, the profits of the three parties fluctuate and lack more severely than the other two models.

To sum up, when the manufacturer is the market leader, the system is also relatively stable. The profit fluctuation and loss of all parties are relatively small. Regardless of the model, when low-carbon sensitivity and channel preference increase too much, the system will fall into chaos. In the chaotic region, the average profit of the supply chain players fluctuates sharply, and it is harmful to market operation.

\section{Chaos control}

Intense market competition may cause uncertainty of the results, have a certain impact on the market pricing system and profits, and even lead to the withdrawal of market players. Chaos always damages the stability of the supply chain and has a negative impact on all competitors. It is harmful to the whole system and economic market. Therefore, the chaos control of the supply chain model has aroused wide attention in academia. In this study, two kinds of control ways are adopted to regulate the disorder of the system.

\subsection{Parameter adjustment control method}

Managers can control chaotic behavior by introducing control parameters $\upsilon$. We set the control parameter to be $\upsilon=0.6$ with the same values of other parameters as the original settings. The parameter adjustment control method is put forward. The system is obtained as follows.

\begin{equation}
\left\{\begin{array}{l}
k_{1, t+1}=(1-v)\left(k_{1, t}+g_{1} k_{1, t} \frac{\partial \pi_{r 1, t}}{\partial k_{1, t}}\right)+v k_{1, t} \\
k_{2, t+1}=(1-v)\left(k_{2, t}+g_{2} k_{2, t} \frac{\partial \pi_{r 2, t}}{\partial k_{2, t}}\right)+v k_{2, t} \\
w_{1, t+1}=(1-v)\left(w_{1, t}+g_{3} w_{1, t} \frac{\partial \pi_{m, t}}{\partial w_{1, t}}\right)+v w_{1, t} \\
w_{2, t+1}=(1-v)\left(w_{2, t}+g_{4} w_{2, t} \frac{\partial \pi_{m, t}}{\partial w_{2, t}}\right)+v w_{2, t}
\end{array}\right.
\end{equation}

\subsection{Delayed feedback control method}

When $g_1=4$, the system is in a state of chaos, and the increase in adjustment speed of retailer 1's sales commission will make the system fall into chaos. According to the delay control method, the partial information of the output signal after a time delay is fed back to the system as an external input. The delayed feedback control system in this work can be shown in Eq(30):

\begin{equation}
\left\{\begin{array}{l}
k_{1, t+1}=k_{1, t}+g_{1} k_{1, t} \frac{\partial \pi_{r 1, t}}{\partial k_{1, t}}+Z\left[k_{1, t}-k_{1, t+1}\right] \\
k_{2, t+1}=k_{2, t}+g_{2} k_{2, t} \frac{\partial \pi_{r 2, t}}{\partial k_{2, t}}+Z\left[k_{2, t}-k_{2, t+1}\right] \\
w_{1, t+1}=w_{1, t}+g_{3} w_{1, t} \frac{\partial \pi_{m, t}}{\partial w_{1, t}}+Z\left[w_{1, t}-w_{1, t+1}\right] \\
w_{2, t+1}=w_{2, t}+g_{4} w_{2, t} \frac{\partial \pi_{m, t}}{\partial w_{2, t}}+Z\left[w_{2, t}-w_{2, t+1}\right]
\end{array}\right.
\end{equation}

Where, $Z$ is the control factor.

\subsection{Numerical analysis of two control methods}

It can be seen from the above analysis that the manufacturer-dominated market is relatively stable. Therefore, this section focuses on the numerical simulation analysis adopting the parametric adjustment and the delayed feedback control method for the model NG and model RS.

%Í¼8 abcd

\begin{figure}[htbp]
\centering

\subfigure[ Parameter adjustment, Model NG, $g_1$=2.9.]{
\begin{minipage}[t]{0.5\linewidth}
\centering
\includegraphics[width=9cm]{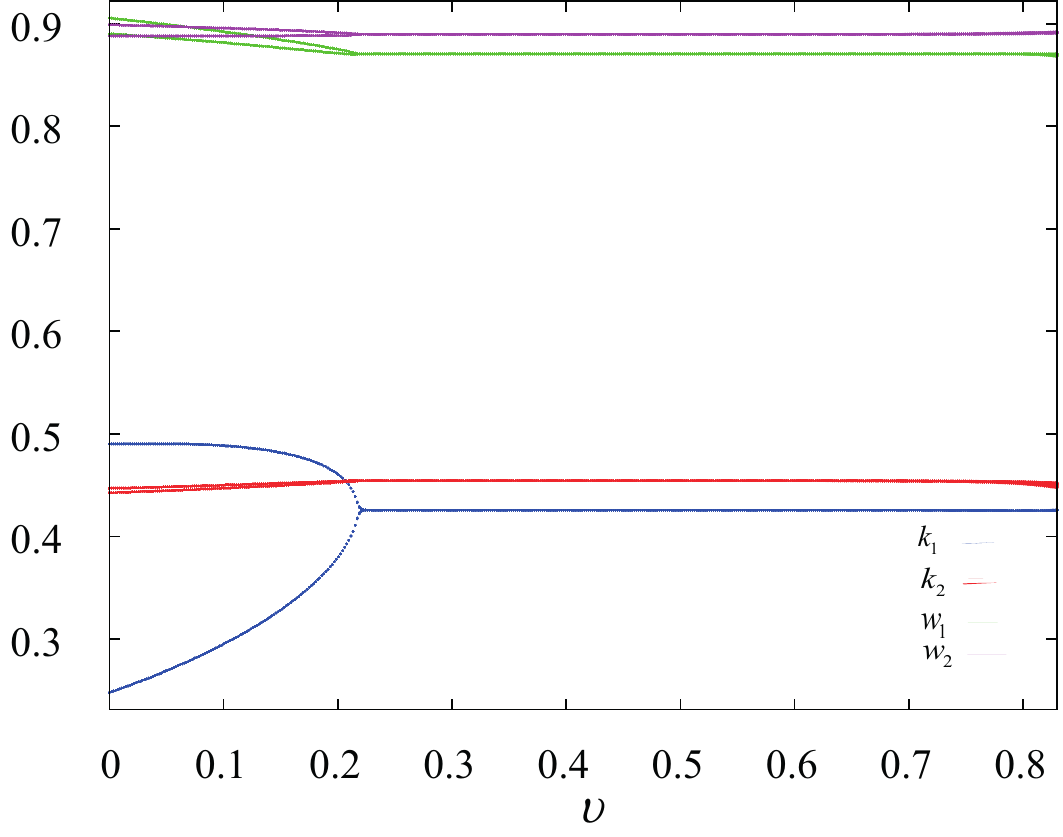}
%\caption{fig1}
\end{minipage}%
}%
\subfigure[ Delayed feedback, Model NG, $g_1$=2.9.]{
\begin{minipage}[t]{0.5\linewidth}
\centering
\includegraphics[width=9cm]{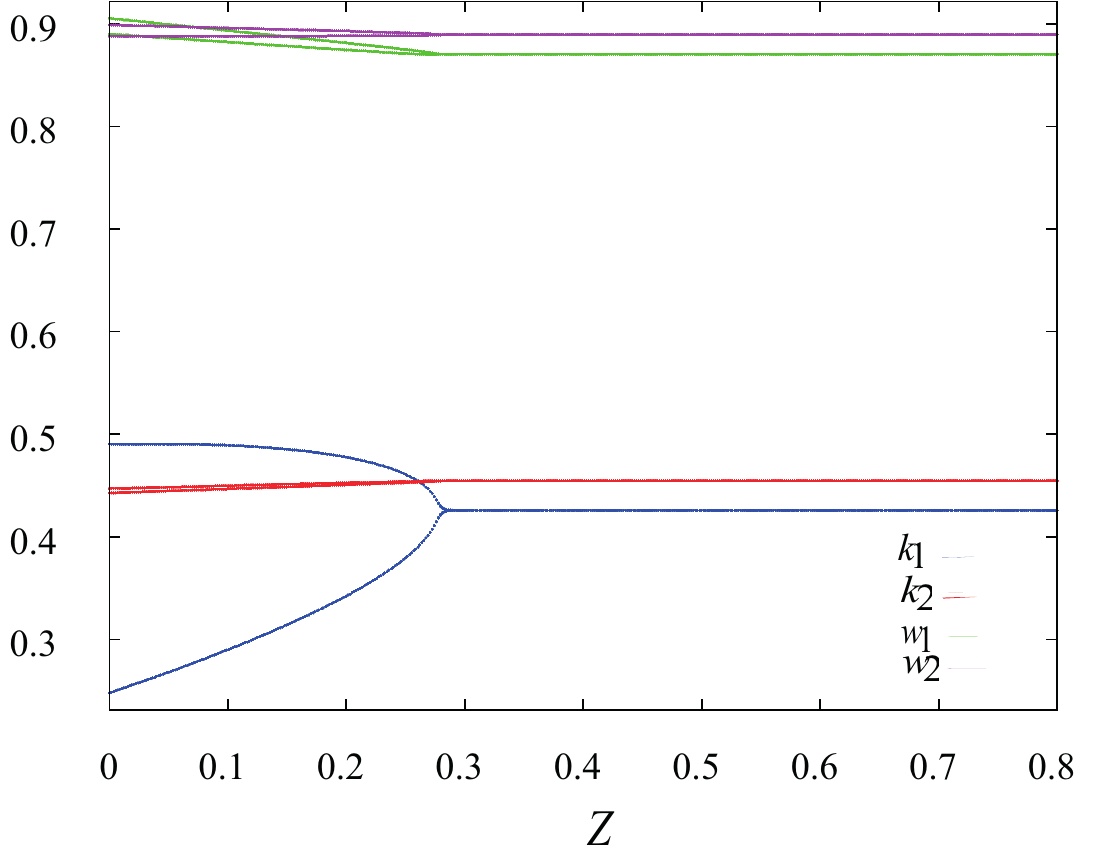}
%\caption{fig2}
\end{minipage}%
}%

\subfigure[ Parameter adjustment, Model RS, $g_1$=4.]{
\begin{minipage}[c]{0.5\linewidth}
\centering
\includegraphics[width=9cm]{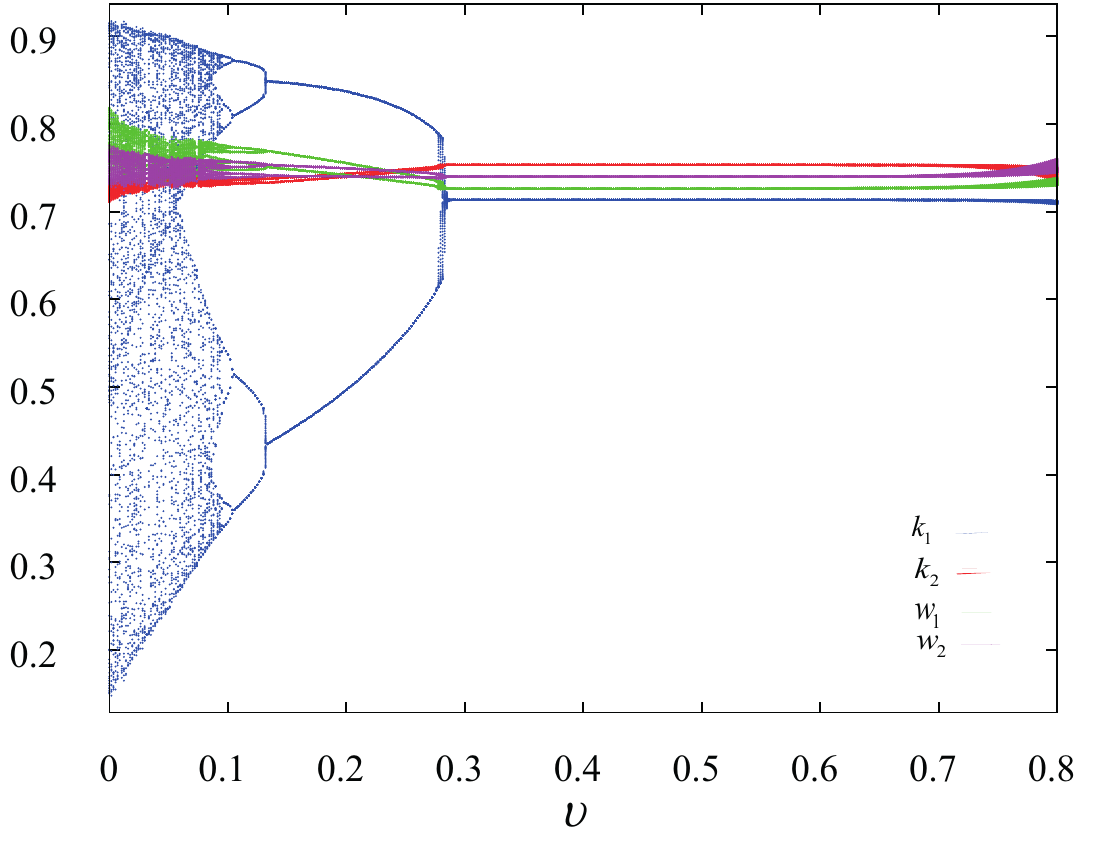}
%\caption{fig2}
\end{minipage}%
}%
\subfigure[ Delayed feedback, Model NG, $g_1$=4.]{
\begin{minipage}[c]{0.5\linewidth}
\centering
\includegraphics[width=9cm]{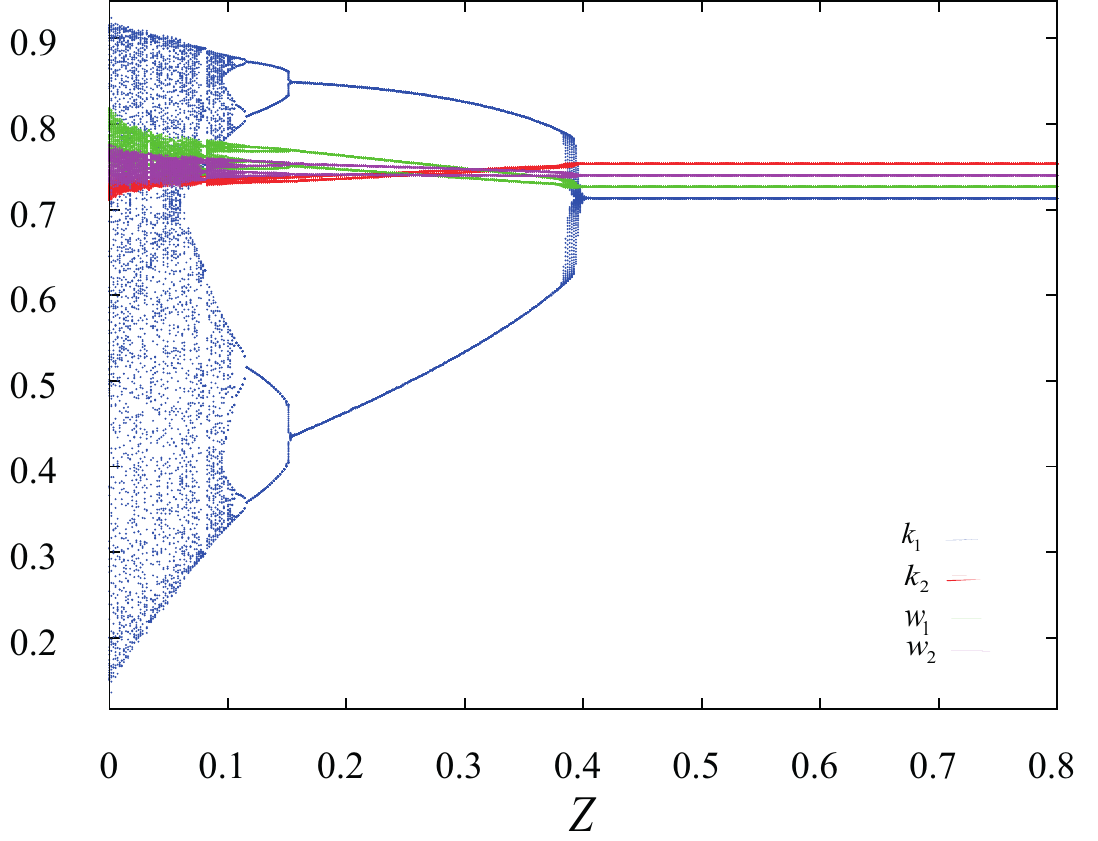}
%\caption{fig2}
\end{minipage}%
}%
\centering
\caption{Bifurcation diagram under two control methods after chaos control.}
\end{figure}

Figure 8 illustrates the bifurcation diagrams after chaos control. Figure 8(a) shows the bifurcation diagram after chaos control using the parameter control method in model NG. When $g_1=2.9$ and the other parameters remain unchanged, the system is in an unstable state, but after introducing the control parameter $\upsilon$, the period-doubling bifurcation state of the system is controlled. When $\upsilon=0.22$, the system will enter into the stable state. In addition, Figure 8(b) illustrates the bifurcation diagram after adopting the delayed feedback control method in model NG. With the increase of $Z$, the system gradually returns to the stability region. And it finally reaches equilibrium at $Z=0.27$.  When $g_1=4$, Figure 8(c) shows that the system is in the chaotic state in model RS. The system is gradually controlled using the parameter control method. When $Z=0.29$, the system returns to the stable state. Figure. 8 (d) presents that the system is controlled at $g_1=4$ in the model RS, when the delayed feedback control method is adopted. The system is restored to the stable state when $Z=0.4$. It is shown that both control methods can control chaos well.

%Í¼Æ¬9 abcd

\begin{figure}[htbp]
\centering

\subfigure[ Model NG, before chaos control.]{
\begin{minipage}[t]{0.5\linewidth}
\centering
\includegraphics[width=9cm]{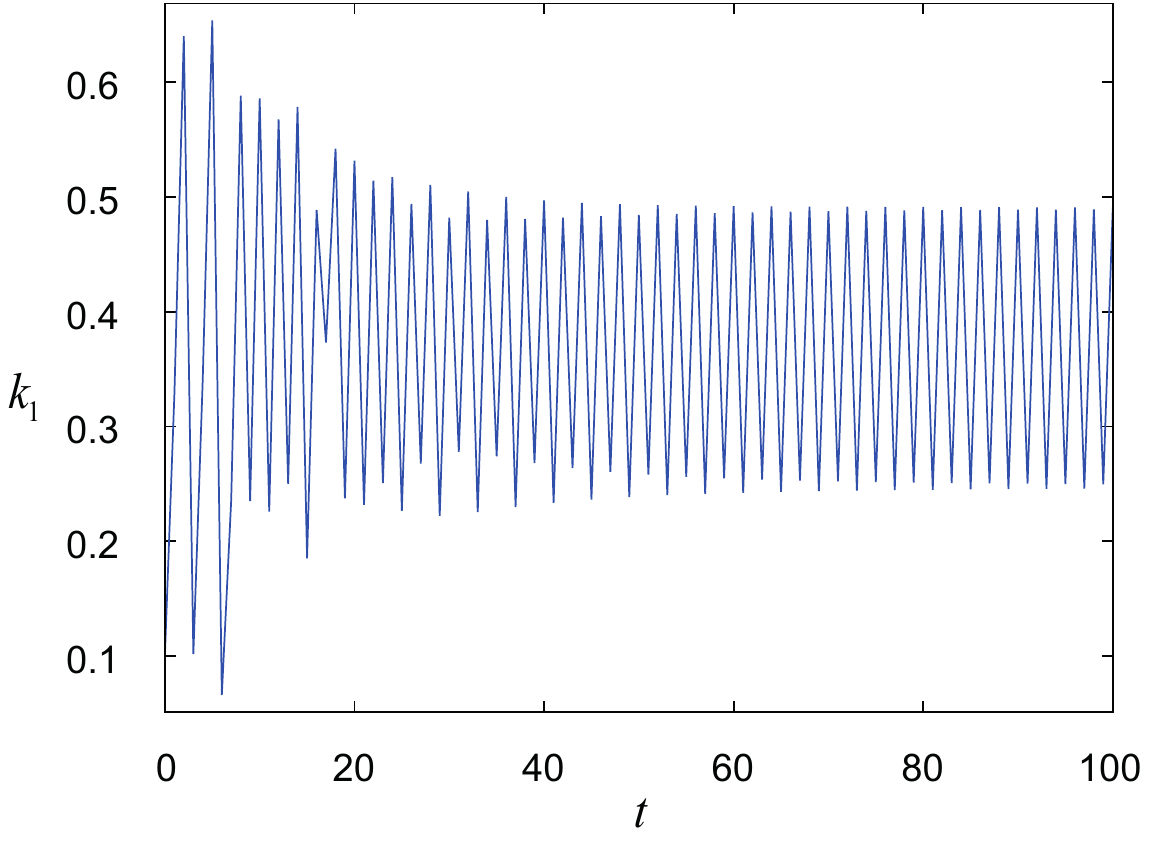}
%\caption{fig1}
\end{minipage}%
}%
\subfigure[ Model NG, after chaos control.]{
\begin{minipage}[t]{0.5\linewidth}
\centering
\includegraphics[width=9cm]{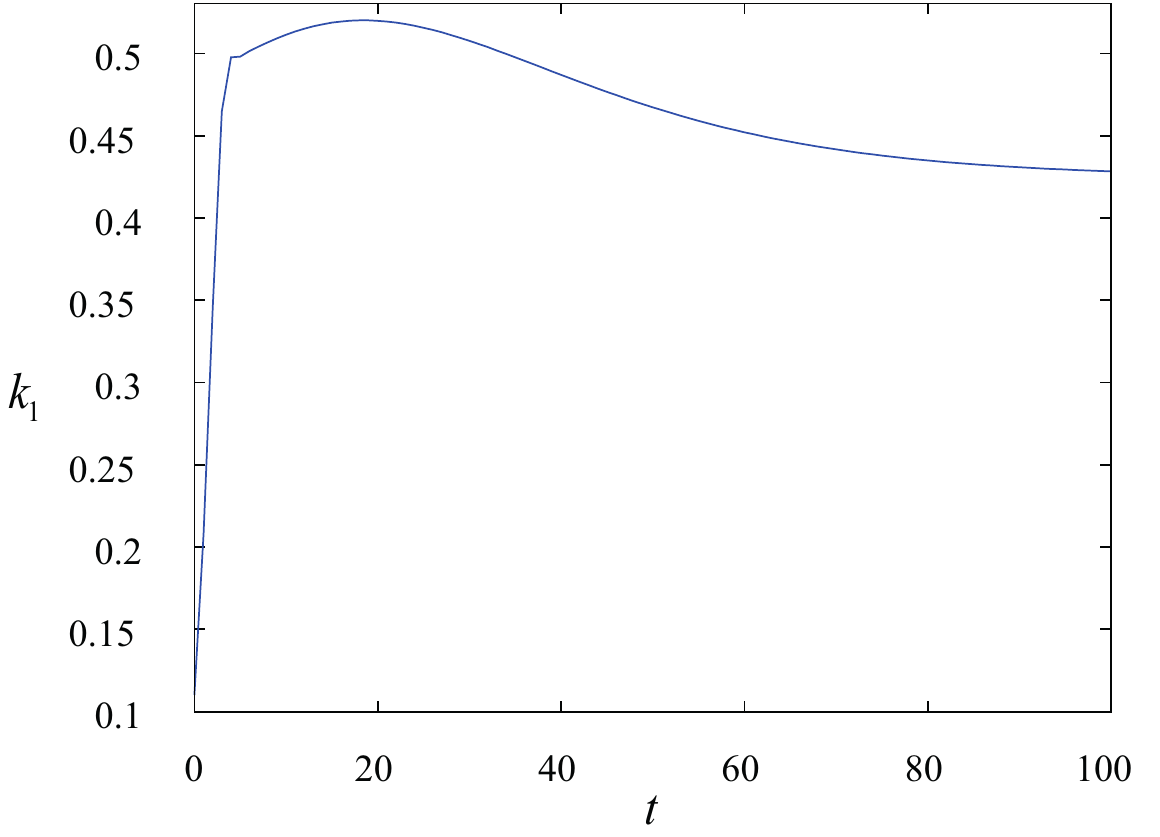}
%\caption{fig2}
\end{minipage}%
}%

\subfigure[ Model RS, before chaos control.]{
\begin{minipage}[c]{0.5\linewidth}
\centering
\includegraphics[width=9cm]{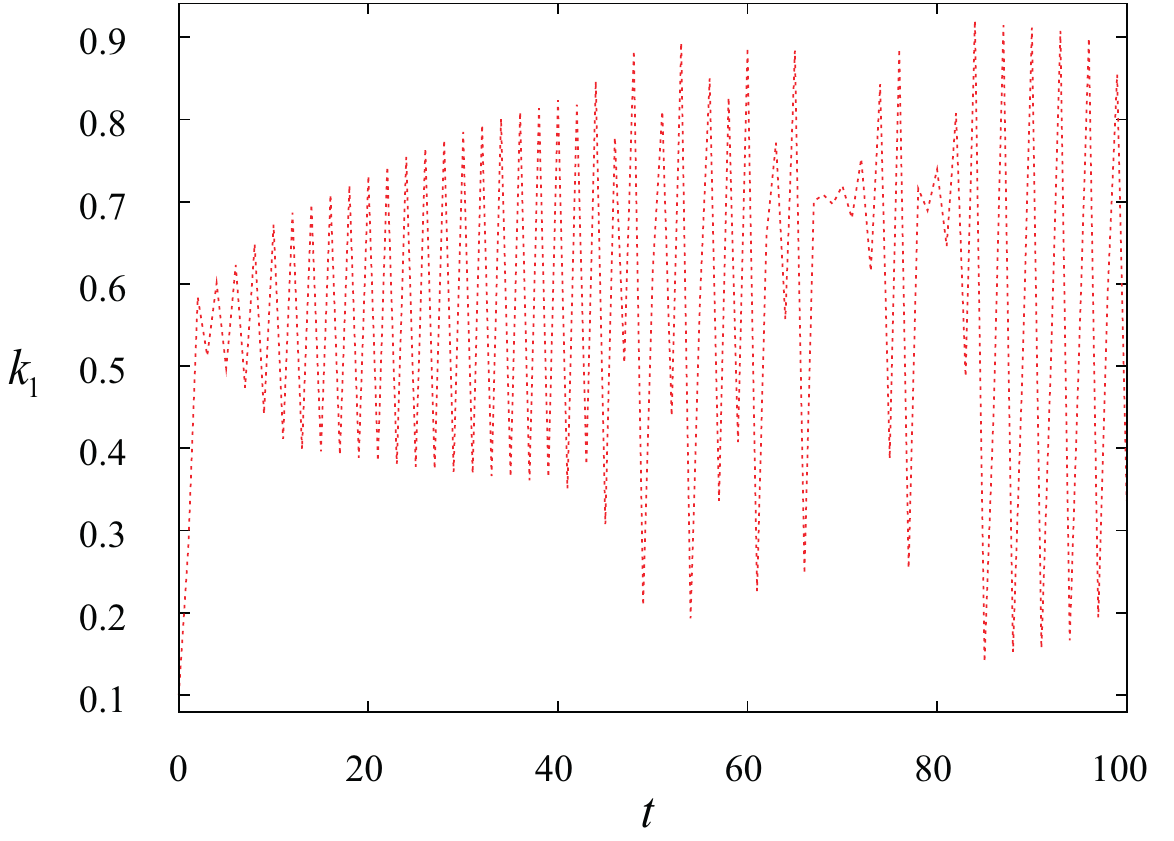}
%\caption{fig2}
\end{minipage}%
}%
\subfigure[ Model RS, after chaos control.]{
\begin{minipage}[c]{0.5\linewidth}
\centering
\includegraphics[width=9cm]{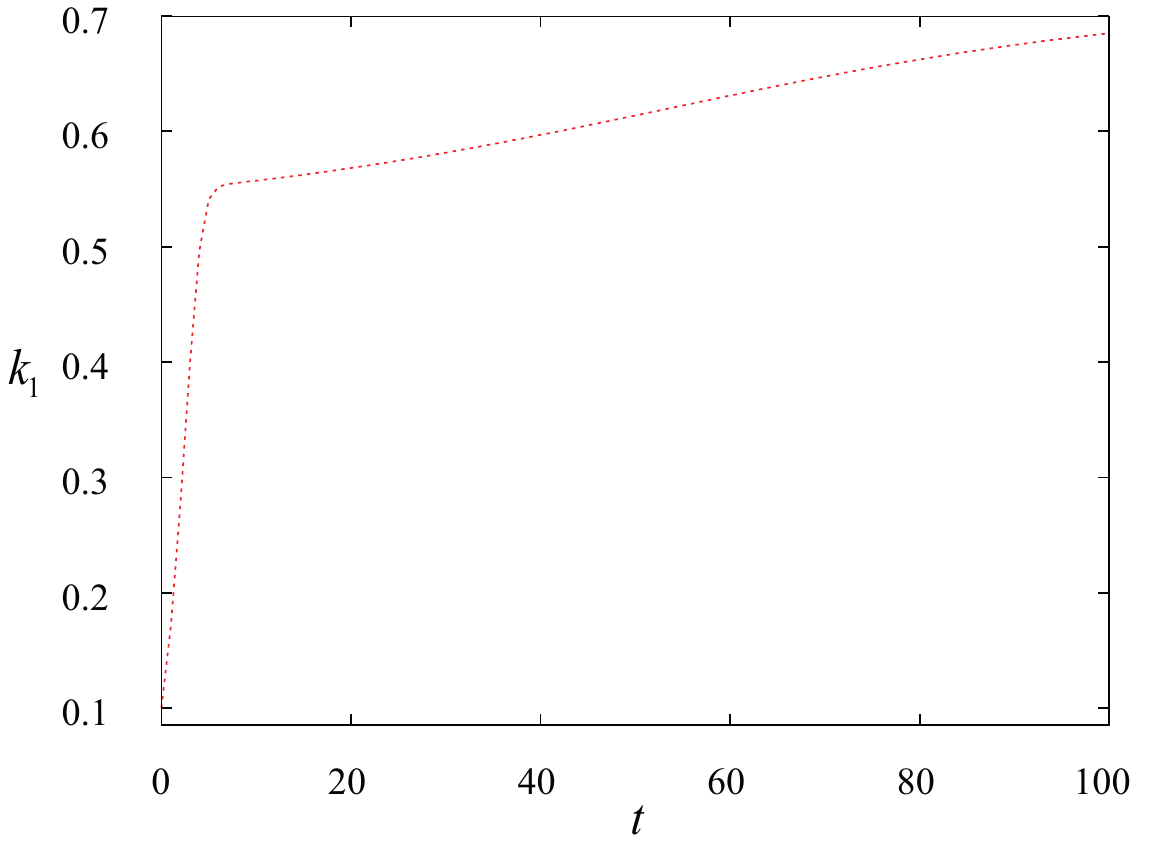}
%\caption{fig2}
\end{minipage}%
}%
\centering
\caption{Time series diagrams after chaos control.}
\end{figure}

Figure 9 presents the time series diagrams of model NG and model RS after chaos control using the parameter control method. As shown in Figure 9 (a), in model NG, when $g_1=2.9$, the system suffers a drastic fluctuation in the chaotic state. Therefore, during the period from 0 to 100, it fluctuates violently. Figure 9 (b) is plotted to explain the time series diagram after the system is controlled in model NG. It can be seen that $k_1$ no longer fluctuates drastically with period $t$, and the system is well controlled. Figure 9(c) shows that in model RS, when $g_1=4$, the system undergoes chaos. The same violent fluctuations occur as $k_1$ varies from 0 to 100 period. Figure 9(d) shows the situation after the system is controlled. It can be seen from Figure 9 that in model NG and RS, the parameter adjustment control method help the system restore stability.

\section{Conclusion}

In the presented paper, we set up a dual-channel supply chain model including a manufacturer, an online retailer, and an offline retailer. And we further investigate complex characteristics in the pricing game process for the supply chain in which manufacturer provides retailers with the same low-carbon products under the government's subsidy.

The game behaviors of the supply chain system under three power structures in a single period are compared and analyzed. It is easy to verify that the growth in two types of consumers' preferences has positive effects on the wholesale price and the sales commission in the three cases.

Through numerical simulation, we further discuss the complex characteristics of the dual-channel supply chain in the dynamic game process under three power structures. The results show that when the subsidy is implemented by the government and the manufacturer becomes the market leader, the price system is relatively more stable and profits are greater. Under the three models, the competition between two retailers will produce substitutability between the two sales channels and make the manufacturer profitable. But when the market is in a highly competitive situation, traditional channel and online channel can mutually substitute each other. Then, the system of the dual-channel supply chain will also fall into chaos.

If consumers prefer online sales channel, wholesale prices and sales commissions will increase. Too strong preference will give rise to instability and even chaos. In this situation, when the manufacturer with the government low-carbon subsidy dominates the market, the system is relatively stable. Consumers have more environmental awareness and prefer low-carbon commodities, which is contribute to manufacturer and retailers increasing wholesale prices and sales commissions. However, when the players' decision-making simultaneously or retailers dominate the market, if consumers excessively pursue the low-carbon effect, the system is more likely to fall into bifurcation and chaos. Chaos will lead to the loss of profits for both manufacturer and retailers, resulting in market disruption. The system can be successfully controlled by utilizing two different approaches to control. With the change of control parameters and factors, the system returns to stability.

Regarding future research directions, the impact of some related random factors on the supply chain channels' equilibrium strategies can be further studied. In addition, the influences of different carbon policies in the multi-channel supply chain system also can be considered in future research. Last but not least, different types of contract and inventory strategies can also be considered in the model for comparative study.

\section*{Acknowledgement(s)}

The authors would like to thank the reviewers for their careful reading and some pertinent suggestions. The research was supported by the project of Tianjin University of Technology and Education (KYQD202227), and supported by Tianjin Educational Science Planning Project (CAE220038).

\section*{Disclosure statement}

Declarations of interest: none.

\bibliographystyle{ws-ijbc}
%\biboptions{authoryear}
\bibliographystyle{plain}
\bibliography{reference}

%\end{multicols}
\end{document}